\pdfoutput=1
\documentclass[]{article}
\usepackage{epstopdf}% To incorporate .eps illustrations using PDFLaTeX, etc.
\usepackage{subfigure}% Support for small, `sub' figures and tables

\usepackage[numbers,sort&compress]{natbib}% Citation support using natbib.sty
\RequirePackage[OT1]{fontenc}
\usepackage{amssymb,makeidx,latexsym,epsf,amssymb,graphicx, url, setspace, amsthm,amsmath,booktabs}

\begin{document}

\begin{titlepage}
\noindent
{\Large\bf  On multivariate  modifications of Cramer Lundberg risk model with constant intensities}\\

 \vspace{1cm}
 \noindent
\uppercase{{\Large Pavlina K.~Jordanova and Milan Stehl\'\i k}}\\

\vspace{1cm}
\noindent
{\it   Faculty of Mathematics and Informatics, Shumen University, Shumen, Bulgaria\\
and\\
Department of Applied Statistics and Linz Institute of Technology, Johannes Kepler University in Linz, Austria\\
Department of Statistics, University of Valparaiso, Valpara\'iso,  Chile\\
}

\noindent
{Address correspondence to Milan Stehl\'\i k,  ........ ~~Email:  milan.stehlik@jku.at, ~~Tel.+56 32 2654680, +4373224686808; Fax: +4373224686800}
\end{titlepage}

\noindent
{\Large\bf On multivariate  modifications of Cramer Lundberg risk model with constant intensities}

\begin{abstract}
The paper considers very general multivariate modifications of Cramer-Lundberg risk model. The claims can be of different types and can arrive in groups. The groups arrival processes within a type have constant intensities. The counting groups processes are dependent multivariate  compound Poisson processes of type I. We allow empty groups and show that in that case we can find stochastically equivalent Cramer-Lundberg model with non-empty groups.

The investigated model generalizes the risk model with common shocks, the Poisson risk process of order k, the Poisson negative binomial, the Polya-Aeppli, the Polya-Aeppli of order k among others. All of them with one or more types of polices.

The relations between the numerical characteristics and distributions of the components of the risk processes are proven to be corollaries of the corresponding formulae of the Cramer-Lundberg risk model.
\end{abstract}

{\bf Keywords:}  {Cramer-Lundberg risk model, Multivariate risk processes, risk process approximations.}\\

{\bf Subject classification codes:}  primary 62P05, secondary 60G10.

\section{Introduction}
A basic model in collective insurance risk theory is the one, introduced by Filip Lundberg \cite{lundberg1932some} and Harald Cramer
\cite{cramer1930mathematical}. It is called Cramer-Lundberg risk model. Google search engine finds approximately 173 000 results in 0,45 seconds,
which speaks about the popularity of this model. In the next section we remind the well known facts about this model that we use further on.
 They can be found in many textbooks  in risk theory, e.g. Grandell \cite{grandell2012aspects}, Embrechts, Kl\"uppelberg, Mikosch \cite{embrechts2013modelling})  Rolski, Schmidli, Schmidt,   Teugels \cite{rolski2009stochastic}, Gerber \cite{gerber1979introduction}.
However, addressing of the possibly empty groups is not so well developed. In order to overcome this gap, we consider possibly empty claims
in Section 3. It can appear for example if a claim arrive but the insurer decide that there is no enough reasons to pay for it.
We apply the
previously mentioned results in the case when the claims can be of different types and can arrive in groups(never mind possible empty or not)
and the inter-arrival times between groups are exponential and independent identically distributed(i.i.d.). This model is discussed in Section 4.

 In the end of the paper we show that the Poisson model with common shock introduced in Cossette and Marceau \cite{A1} is its special case of the model defined here, with a generalization of Marshall-Olkin arrival process \cite{MO1967, MO1988}. We would obtain a stochastically equivalent of their model if no simultaneously empty groups are possible, no more then one claim can form a group and $d = 3$.

Wang and Yuen \cite{A3} partially investigate and generalize the model of Cossette and Marceau \cite{A1} for non-homogeneous group counting processes.

Under the constrain that the claim sizes ar Erlang see Yuen, Guo, and Wu \cite{A2}.

In Sections 5 and 6 some particular cases follow.
First we  show that if the claims of different types arrive in groups and the groups counting processes
are independent Compound Poisson processes with possibly different intensities we obtain particular case of the model defined in Section 4.

In the last section we explain why our novel model of risk process, generalizes the Compound Compound Poisson risk process,
Poisson risk process of order k, Poisson negative binomial, Polya-Aeppli and Polya-Aeppli of order k
 risk models, introduced in series of works of Minkova, Kostadinova and Chukova\cite{chukova2015polya, KrasiDoctoralConference, kostadinova2013poisson, kostadinovapolya, kostadinova2014bivariate, minkova1900polya, minkova2009compound}, the Poisson model with common shock introduced in Cossette and Marceau \cite{A1} among others. All of them with one or more types of polices.
%\textcolor[rgb]{1.00,0.00,0.00}{The corresponding probabilities of ruin in finite time and for different, but fixed initial capitals, are simulated
%for different but fixed distribution of the claim sizes with finite variance  at the ends of   the  sections.}
%The last approach is used e.g. in Jordanova(2005)\cite{} and Mitov(2009)\cite{} .

Through the paper $l_{X}$ denotes the Laplace - Stieltjes transform of the random variable (r.v.) or the function $X$ and $g_{X}$ is the probability generating function(p.g.f.) of the r.v. $X$. $F^{n*}(x)$ is the $n$-th convolution of $F$, and by convention $F^{0*}(x) = I\{x \geq 0\}$, $x \in \mathbb{R}.$

\section{The Cramer-Lundberg risk model - revisited}
The Cramer-Lundberg risk model is usually described via the following conditions (C1) - (C5). See Figure \ref{fig:Fig1}.

  \begin{center}
\includegraphics[width=1\textwidth]{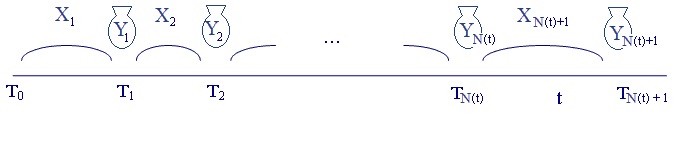}

{\small Fig. 1. The claim arrival process. \label{fig:Fig1}}
\end{center}

\begin{itemize}
  \item [(C1)] The clime sizes $Y_1, Y_2, ...$ are positive independent identically distributed (i.i.d.) random variables (r.vs.) having common non-lattice cumulative distribution function (c.d.f.), with mean $\mathbb{E}(Y_1) < \infty$ and variance $\mathbb{D}(Y_1) \leq \infty$.
  \item [(C2)] The claims occur at random time points $0  < T_1 < T_2 < ... \quad a.s.$ and by convention $T_0 = 0$.
  \item [(C3)] The claims arrival process is defined as
  $$N(t) = sup\{k \in \mathbb{N}: T_k \leq t\}, \,\, t > 0, \,\, sup\, \emptyset \,\, = 0.$$
  \item [(C4)] The inter-arrival times $X_k = T_k - T_{k-1}$, $k = 1, 2, ..$ are i.i.d. exponentially distributed with $\mathbb{E} X_1 = \frac{1}{\lambda}$, $\lambda > 0$.
  \item [(C5)]The sequences $X_1, X_2, ...$ and $Y_1, Y_2, ...$ are independent.
  \end{itemize}
  It is well known that condition (C4) is equivalent to the condition (C6)
  \begin{itemize}
  \item [(C6)]$N = \{N(t): t \geq 0\}$ is a homogeneous Poisson process(HPP) with intensity $\lambda > 0$.
  Briefly $N \sim HPP(\lambda)$.
   \end{itemize}

   Usually the following two stochastic processes are related with the above model.
\begin{itemize}
  \item [$S$] - the total claim amount process. More precisely it is the sum of the claims that have occurred up to time t $$S(t) = I\{N(t) > 0\} \sum_{i=1}^{N(t)} Y_i, \,\, t \geq 0.$$
  \item [$R$] - the risk  process $R(t) = u + ct - S(t)$, $t \geq 0$,
where $c > 0$ is the premium income rate  and $u \geq 0$ is the initial capital.
\end{itemize}
From practical point of view, it is important to guarantee that the mean income is bigger that the mean expenditures. Therefore the following characteristic is a basic characteristics of the risk models
  \begin{equation}\label{safetiloadingCLM}
\rho = \lim_{t \to \infty} \frac{ER(t)}{ES(t)} = \lim_{t \to \infty} \frac{u + ct - \mathbb{E}S(t)}{\mathbb{E}S(t)} = \frac{c}{\lambda \mathbb{E}(Y_1)} - 1.
\end{equation}
It is called safety loading.

From insurer's point of view, it is interesting to know the description of the time of ruin. It  depends on the initial capital $u$ and  is defined by
$$\tau(u) = inf\{t > 0: R(t) < 0| R(0) = u\}, \,\, u \geq 0, \,\, inf \emptyset\,\, = \infty.$$
The probability for ruin in "infinite horizon", also depends on $u$ and
$$\psi(u) = P(\tau(u) < \infty| R(0) = u).$$
The corresponding survival probability is $\delta(u) = 1 - \psi(u)$.

It is well known that the risk theory is closely related with the theory of random walks. It can be shown that for $Z_0 = 0$ and $Z_k = \sum_{i = 1}^k(Y_i - cX_i)$, $k = 1, 2, ...$
\begin{equation}\label{randomwalk}
\psi(u) = \mathbb{P}(\sup_{i \in \mathbb{N}} Z_i > u).
\end{equation}
Therefore in order to obtain $\lim_{u \to \infty} \psi(u) = 0$  we need to impose  the net profit condition: $\rho > 0$, which is equivalent to the condition $\mathbb{E}(Y_1 - cX_1) < 0$, and to
\begin{equation}\label{NPC}
c > \lambda \mathbb{E}(Y_1).
\end{equation}
Otherwise for any initial capital $u \geq 0$, $\psi(u) = 1$. Because  here we consider "infinite horizon" everywhere we suppose that this condition is satisfied. In that case it is well known that

\begin{equation}\label{0}
\psi(0) = \frac{1}{1 + \rho} = \frac{\lambda \mathbb{E}(Y_1)}{c}, \quad \delta(0) = \frac{\rho}{1 + \rho}.
\end{equation}
{Cor.  11.3.1., p. 471 in the book of Rolski, Schmidli, Schmidt, Teugels(2009) \cite{rolski2009stochastic} give us formula for $\mathbb{E}(\tau(0)|\tau(0) < \infty)$ in the case when $\mathbb{D}Y_1 < \infty$
\begin{equation}\label{tau0}
\mathbb{E}(\tau(0)|\tau(0) < \infty) = \frac{1}{2\lambda \rho}\left(\frac{\mathbb{D}Y_1 }{(\mathbb{E}Y_1 )^2} + 1\right) = \frac{\mathbb{D}Y_1 + (\mathbb{E}Y_1 )^2}{2\mathbb{E}Y_1(c-\lambda \mathbb{E}Y_1)}.
\end{equation}
%and for all $x > 0$
 %\begin{equation}\label{tau0measure}
%\mathbb{P}(x < \tau(0) < \infty) \leq \frac{\lambda[\mathbb{D}Y_1 + (\mathbb{E}Y_1 )^2]}{2x(c-\lambda \mathbb{E}Y_1)^2},
%\end{equation}
}
The theory about these risk models can be developed from  the following integro differential equation:
\begin{equation}\label{deltaPrim}
\delta'(u) = \frac{\delta(u)}{(1+\rho)\mathbb{E}Y_1} - \frac{1}{(1+\rho)\mathbb{E}Y_1} \int_0^u \delta(u - y) dF_{Y_1}(y).
\end{equation}
The solution of { (\ref{deltaPrim})} satisfy the defective renewal equation
\begin{equation} \label{renewaldelta}
 \delta(u) = (1 - \psi(0)) + \psi(0)\int_0^u \delta(u-y) dF_I(y).
\end{equation}
where $F_I(x) := \frac{1}{\mathbb{E}(Y_1)}\int_0^x (1 - F(y)) dy$ is the integrated tail c.d.f.(or the c.d.f. of the equilibrium distribution of $Y_1$) and $\bar{F}_I(u) := 1 - F_I(u)$. Therefore
\begin{equation} \label{renewalpsi}
 \psi(u) = \psi(0)\bar{F}_I(u) - \psi(0)\int_0^u \psi(u-y) d\bar{F}_I(y).
\end{equation}
and the solutions of {(\ref{renewalpsi})} and {(\ref{renewaldelta})} can be presented via the following Beekman’s  convolution series \cite{beekman1968collective} or Pollaczek - Khinchin \cite{Pollaczek, Khinchin} formula, which claims
\begin{eqnarray}
\psi(u) &=&   \left(1 - \psi(0)\right) \sum_{i = 1}^\infty \psi^i(0) [1 - F_I^{i*}(u)], \,\, u \geq 0. \\\label{PH}
\nonumber\delta(u) &=& \left(1 - \psi(0)\right) \sum_{i = 0}^\infty \psi^i(0) F_I^{i*}(u),  \,\, u \geq 0.
\end{eqnarray}
The above formulas (\ref{PH}) will be true for all multivariate models discussed in the paper.
They express the fact that, if we consider a sequence of independent repetitions of an experiment, one trial is running a stochastically equivalent random walk $\tilde{Z}_k$ of $Z_k$ and "success" is "$\sup_{i \in \mathbb{N}} \tilde{Z}_i = 0$", then the $\sup_{i \in \mathbb{N}} Z_i$ coincides in distribution with sum of the final values in the sequence, before the first "success" happen. The probability for success in one trial is $P(\sup_{i \in \mathbb{N}} \tilde{Z}_i = 0) = \delta(0).$
Therefore  $\sup_{i \in \mathbb{N}} Z_i$ is a compound Geometrically distributed, with parameter $\delta(0)$. Moreover if we denote the first upper record time of $Z_1, Z_2, ...$ by $L(1)$. The distribution of the summands coincides with the value of $Z_{L(1)}$ given "$Z_{L(1)} > 0$", which is the distribution of the deficit at the time of ruin with initial capital $0$. More precisely
 $$(Z_{L(1)}|Z_{L(1)} > 0) \stackrel{\rm d}{=} (-R(\tau(0)+)| \tau(0) < \infty)$$ and $$\mathbb{P}(-R(\tau(0)+) \leq x| \tau(0) < \infty) = F_I(x).$$ Then
 \begin{equation}\label{0}
 \mathbb{E}(-R(\tau(0)+)| \tau(0) < \infty) = \frac{\mathbb{E}Y_1^2}{2\mathbb{E}Y_1} = \mathbb{E}(\tau(0)|\tau(0) < \infty)(c-\lambda \mathbb{E}Y_1).
 \end{equation}
  The joint distribution of the severity of (deficit at) ruin and the risk surplus just before the ruin with initial capital zero is:
\begin{equation}\label{jointdeficite}
\mathbb{P}(-R(\tau(0)+) > x, R(\tau(0)-) > y|\tau(0) < \infty) = \bar{F}_I(x + y).
\end{equation}
The distribution of the claim causing ruin is:
\begin{equation}\label{theclaimcausingruin}
\mathbb{P}(R(\tau(0)-)-R(\tau(0)+) \leq x| \tau(0) < \infty) = \frac{1}{\mathbb{E}(Y_1)}\int_0^x y dF_{Y_1}(y).
\end{equation}
See e.g. Embrechts, Kl\"uppelberg, Mikosch \cite{embrechts2013modelling}.

The above considerations show that for initial capital zero we have explicit formulae for the numerical characteristics if the Cramer-Lundberg model, however when we consider arbitrary strictly positive initial capital $u$, the situation is not so simple. in that case, however we have explicit form of the Laplace - Stieltjes transform of $\delta(u)$. It has the form
\begin{equation}\label{LStTrdeltaGeneral}
l_{\delta}(s) = \frac{\frac{\rho}{1 + \rho}}{1 - \frac{1}{1+\rho} \left(\frac{1 - l_{Y_1}(s)}{s\mathbb{E}Y_1}\right)}.
\end{equation}

Moreover for all $u \geq 0$, if $G(u, y) = \mathbb{P}(-R(\tau(u)+) \leq y,  \tau(u) < \infty),$ then
\begin{equation}\label{G0y}
G(0, y) = F_I(y)\psi(0) = \frac{\lambda}{c}\int_0^y \bar{F}_{Y_1}(z) dz,
\end{equation}
\begin{equation}\label{dG}
\frac{\partial}{\partial u} G(u, y) = \frac{\lambda}{c}\left[G(u, y) - \int_0^u G(u - x, y)dF_{Y_1}(x) - [F_{Y_1}(u + y) - F_{Y_1}(u)]\right],
\end{equation}
 and
\begin{equation}\label{G}
G(u, y) = \frac{\lambda}{c}\left[\int_0^u G(u - x, y)[1 - F_{Y_1}(x)]dx + \int_u^{u + y}[1 - F_{Y_1}(x)]dx\right], \,\, u \geq 0.
\end{equation}
See Gerber and Shiu \cite{GSh} or Klugman, Panjer  and  Willmot \cite{klugman2012loss}.

%For $r > 0$, such that $E e^{r Y_1} < \infty$, denote by $g(r) = \lambda[E e^{r Y_1} - 1] - cr$. Then, $g(0) = 0$, $g'(0) < 0$, $g''(r)> 0$, $Ee^{-r(R(t) - u)} = e^{tg(r)}$, the process $\{e^{-r[R(t)-u]-tg(r)}: t \geq 0\}$ is a martingale and for all $T > 0$ and $r > 0$,
%\begin{eqnarray}
%\nonumber \psi(u, T) &\leq& \frac{e^{-ru}}{\mathbb{E}[e^{-r [R(\tau(u)) - u] - g(r)\tau(u)}| \tau(u) < T]} \\
%&\leq& \label{psiT} \frac{e^{-ru}}{\mathbb{E}[e^{-g(r)\tau(u)}| \tau(u) < T]} \leq e^{-ru}\sup_{y \in [0, T]}e^{yg(r)}.
%\end{eqnarray}

A relatively good approach to $\psi(u)$ for arbitrary $u \geq 0$ is via the so called Lundberg exponent. Given, the small claim Cramer-Lundberg condition is satisfied, or more precisely, if these exists the Cramer-Lundberg exponent $\epsilon > 0$:
\begin{equation}\label{Lundbergcondition}
\mathbb{E} e^{\epsilon(Y_1 - cX_1)} = 1,
\end{equation}
 or, which is the same,
 \begin{equation}\label{LE}
 \int_0^\infty e^{\epsilon x} d  F_I(x) = 1 + \rho,
 \end{equation}
  then
\begin{equation}\label{psiLundbergexponent}
\psi(u) = \frac{e^{-\epsilon u}}{\mathbb{E}[e^{-\epsilon [R(\tau(u))-u]}| \tau(u) < \infty]} \leq e^{-\epsilon u}.
\end{equation}

The inequality (\ref{psiLundbergexponent}) allows us to chose $\alpha \in (0, 1)$, then to determine $\epsilon$ by the equality $\epsilon = -\frac{log(\alpha)}{u}$. Given the distribution of $Y_1$, by  equation (\ref{LE}) we can determine the safety loading $\rho$. Finally via (\ref{0}) one can obtain the premium income rate $c$ in such a way that $\epsilon$ will be the Lundberg exponent of the considered
model and $\psi(u) \leq \alpha$.
See  Klugman, Panjer  and  Willmot \cite{klugman2012loss}.

 If additionally
 \begin{equation}\label{condapproxpsi}
  \int_0^\infty x e^{\epsilon x}dF_I(x) < \infty,
\end{equation}
then the following Cramer-Lundberg approximation of the probability of ruin is true
\begin{equation}\label{approxpsi}
\lim_{u \to \infty} e^{\epsilon u}\psi(u) = \frac{\rho}{\epsilon \int_0^\infty x e^{\epsilon x} dF_I(x)}.
\end{equation}
See e.g. Embrechts, Kl\"uppelberg, Mikosch \cite{embrechts2013modelling}.

{The condition (\ref{condapproxpsi}) seems to be very restrictive and is not satisfied by the most of the distributions which can be met in practice. Therefore the following result of Goldie and  Klueppelberg \cite{B1} seems to be very useful. They state that if the net profit condition is satisfied, the integrated tail distribution of $Y_1$ is subexponential if and only if $\delta(u)$ is subexponential and in that case we have the following asymptotic of the probability of ruin
\begin{equation}\label{AsymptoticHeavy}
\lim_{u \to \infty} \frac{\psi(u)}{\bar{F}_I(u)} = \frac{1}{\rho}.
\end{equation}
If the claim sizes belong to the distributional class with regularly varying tail, then the random sum that express the total claim amount within a group will be also regularly varying. Embrechts, Klueppelberg Mikosch \cite{embrechts2013modelling} show that in that case the integrated tail distribution of $Y_1$ will be subexponential and we can apply the above theorem. Moreover given the net profit condition is satisfied, they show that it is sufficient that the distribution of the claim sizes  belongs to the class of dominantly varying distributions
$$D = \{F\,\, c.d.f.\,\, on\,\, (0, \infty):\,\, \limsup_{x \to \infty} \frac{\bar{F}(x/2)}{\bar{F}(x)} < \infty\}$$
in order to achieve asymptotic (\ref{AsymptoticHeavy}).}

The case, when the claim sizes are  exponentially distributed with mean $\mu$ is very well investigated. Then the equation (\ref{Lundbergcondition}) is solved and
\begin{equation}\label{expLundbergexponent}
 \epsilon = \mu( 1 - \psi(0)) = \frac{\rho \mu}{1 + \rho}.
\end{equation}
The solution of (\ref{renewalpsi}) and explicit form of $\psi(u)$ in (\ref{PH}) is
\begin{equation}\label{expPsi}
 \psi(u) = \psi(0) e^{-\epsilon u} = \frac{1}{1 + \rho} e^{-\epsilon u}, \,\, u \geq 0.
\end{equation}
In that case it is easy to see that the solution of (\ref{renewaldelta}) and explicit form of $\delta(u)$ in (\ref{PH}) is
\begin{equation}\label{expDelta}
 \delta(u) = 1 - \psi(0) e^{-\epsilon u} = 1 -  \frac{1}{1 + \rho} e^{-\epsilon u}, \,\, u \geq 0.
\end{equation}
The last two results seems to be obtained by Cramer \cite{cramer1955collective}.
{ For $u \geq 0$,
\begin{equation}\label{tauuExpclaims}
\mathbb{E}(\tau(u)|\tau(u) < \infty) = \frac{c + \lambda u}{c(c\mu-\lambda)}.
\end{equation}
See Rolski, Schmidli, Schmidt, Teugels \cite{rolski2009stochastic}.}

The discussion on applicability of several approximations to probability of sum of claims can be found in \cite{RePEc:prg:jnlpep:v:2014:y:2014:i:3:id:488:p:349-370}. Risk related to dividends of insurance companies is studied in \cite{Divident}.
\vspace*{12pt}

\section{Considerations on the Cramer-Lundberg risk model with possibly empty claims}
It is easy to generalize the above results for the case when $Y_1$ is non-negative r.v. with $P(Y_1 = 0) = p_0 \geq 0$. In that case we need to replace the parameter $\lambda$ in the HPP in (C6) with $\lambda(1 - p_0)$. The last property of the HPPs is well known as thinning(splitting) property of Poisson processes \cite{DaleyVereJones}. Here and further on f.d.d. means "in the sense of the finite dimensional distributions". Although the results in the following lemma are easy to obtain, they play an important role in our next considerations.

{\bf Lemma 3.1} Assume $\tilde{\lambda} > 0$,  $\eta_1, \eta_2, ...$ are i.i.d. non-negative r.vs. with $\mathbb{P}(\eta_1 = 0) = p_0 \geq 0$ and $\tilde{N} \sim HPP(\tilde{\lambda})$, independent on $\eta_1, \eta_2, ...$.  Let $Y_1, Y_2, ...$ be i.i.d. positive r.vs. with $\mathbb{P}(Y_1 \leq x) = \mathbb{P}(\eta_1 \leq x | \eta_1 > 0)$. Suppose $N \sim HPP[\tilde{\lambda}(1 - p_0)]$,  independent on $Y_1, Y_2, ...$. Then
$$\mathbb{P}(Y_1 \leq x) = \frac{\mathbb{P}(\eta_1 \leq x) - p_0}{1 - p_0}, \quad \mathbb{E} e^{-xY_1} = \frac{\mathbb{E}e^{-x\eta_1} - p_0}{1 - p_0},$$
$$\mathbb{E}Y_1 = \frac{\mathbb{E}\eta_1 }{1 - p_0}, \quad \mathbb{D}Y_1 = \frac{\mathbb{D}\eta_1}{1 - p_0} -\frac{p_0(\mathbb{E}\eta_1)^2}{(1 - p_0)^2}.$$
 and
$$I\{\tilde{N}(\cdot) > 0\}\sum_{i = 1}^{\tilde{N}(\cdot)} \eta_i \stackrel{\rm f.d.d.}{=} I\{N(\cdot) > 0\}\sum_{i = 1}^{N(\cdot)} Y_i.$$

{\it Proof.} The first equality follows from the definition of conditional probability. Let us prove the second one.
$$\mathbb{E} e^{-xY_1} = \int_{0^+}^\infty e^{-xy}d\mathbb{P}(\eta_1 \leq x | \eta_1 > 0) = \int_{0^+}^\infty e^{-xy}d\frac{\mathbb{P}(\eta_1 \leq x) }{1 - p_0} = \frac{\mathbb{E}e^{-x\eta_1} - p_0}{1 - p_0}.$$
In order to compare the distributions, for all $t \geq 0$ we use the uniqueness of the correspondence between the distribution and its Laplace - Stieltjes transform.
\begin{eqnarray*}
  \mathbb{E} e^{-x I\{\tilde{N}(t) > 0\}\sum_{i = 1}^{\tilde{N}(t)} \eta_i}  &=& e^{-\tilde{\lambda}t(1 - \mathbb{E} e^{-x\eta_1})} = e^{-\tilde{\lambda}t(1 - p_0)(\frac{1 - \mathbb{E} e^{-x\eta_1}}{1 - p_0})} =\\
   &=& e^{-\tilde{\lambda}t(1 - p_0)(\frac{1 - \mathbb{E} e^{-x\eta_1} \pm p_0}{1 - p_0})} =\\
   &=& e^{-\tilde{\lambda}t(1 - p_0)(1 -  \frac{\mathbb{E}e^{-x\eta_1} - p_0}{1 - p_0})} =\\
   &=& e^{-\tilde{\lambda}t(1 - p_0)(1 - \mathbb{E} e^{-xY_1})} = \mathbb{E} e^{-x I\{N(t) > 0\}\sum_{i = 1}^{N(t)} Y_i}.
\end{eqnarray*}
Due to the independence and homogeneity of the additive increments of homogeneous Poisson process, the analogous equalities could be proven for the additive increments and consequently for all finite dimensional distributions.

Therefore if the claims in the model $C_1 - C_6$ are possibly empty, we can find a stochastically equivalent Cramer-Lundberg model without empty claims, and applying the well known formulae for that model we obtain the numerical characteristics of the model $C_1 - C_6$ with possibly empty claims. In our considerations instead of only one claim $Y_1$ we will be the total claim amount within a group, and the claims will be possible to have $d \in N$ different types.

\section{The Cramer-Lundberg model with d types of business, and  with multivariate compound Poisson of type I groups arrival process.}

Let us now assume that the claims are of $d \in N$ different types and it is possible they arrive in groups.
Denote by $\overrightarrow{\tilde{Y}}_{ij}$ the $j$-th claim size within $i$-th group, $i, j = 1, 2, ...$. Due to our assumption that $d$ types of claims are possible, $\overrightarrow{\tilde{Y}}_{ij}$ is a random vector,
$\overrightarrow{\tilde{Y}}_{ij} = (\tilde{Y}_{ij1}, \tilde{Y}_{ij2}, ...., \tilde{Y}_{ijd}),$
where
$\tilde{Y}_{ijs}$ represents $j$-th claim size of type $s$, within $i$-th group, $s = 1, 2, ..., d$. By assumption $\mathbb{P}(\overrightarrow{\tilde{Y}}_{ij} = (0, 0, ..., 0)) = 0$. The claim amounts within a group are independent on the inter-arrival times between groups. The claim sizes of a given type are assumed to be i.i.d. Therefore $\tilde{Y}_{ijs} \stackrel{\rm d}{=} \tilde{Y}_{11s}$ for all $i, j \in \mathbb{N}$ and in cases when it is clear we will skip the first index.

$\tilde{N}(t)$ is the number of groups, inclusively possibly empty groups, that have occurred up to time $t > 0$. On Figure \ref{fig:Fig2} the first four lines represent an example of groups arrival processes in the insurance company. Separate diamonds represent individual claims. They are of four different types. The first line is for the groups of claims of the first type. The second line describes the arrivals of the groups of claims of the second type. The third row is for the third type and the forth depicts the arrivals of the claims of the fourth type. The last row represents the merging of these Poisson processes. More theory about transformations of Poisson processes can be found e.g. in \cite{DaleyVereJones}.

\begin{center}
\includegraphics[width=.9\textwidth]{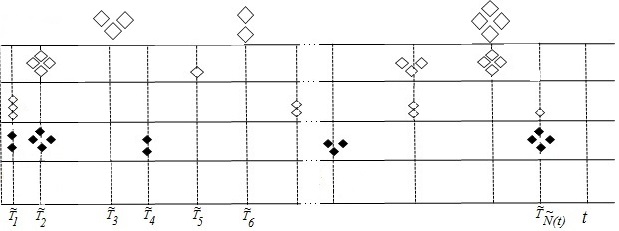}

{\small Fig. 2. Example of a claim arrival processes with Multivariate compound Poisson of type I counting process and $d = 4$.   \label{fig:Fig2}}
\end{center}

We assume that $\tilde{N} = \{\tilde{N}(t): t \geq 0\}$ is a homogeneous Poisson process with intensity $\tilde{\lambda}$.

$\overrightarrow{\tilde{U}}_i$ is the number of claims within $i$-th group. Assume that the random vectors
$\overrightarrow{\tilde{U}}_i = (\tilde{U}_{i1}, \tilde{U}_{i2}, ..., \tilde{U}_{id})$, $i = 1, 2, ...$
 are i.i.d. and independent on $\tilde{N}$. The $s$-th coordinate of this vector represents the number of claims within $i$-th group and of type $s = 1, 2, ..., d$. On the outcome depicted on Figure \ref{fig:Fig2} the vector $\overrightarrow{\tilde{U}}_1$ has coordinates $\tilde{U}_{11} = 0, \tilde{U}_{12} = 0, \tilde{U}_{13} = 3, \tilde{U}_{14} = 2$, $d = 4$. The vector $\overrightarrow{\tilde{U}}_2$ has coordinates $\tilde{U}_{11} = 0, \tilde{U}_{12} = 4, \tilde{U}_{13} = 0, \tilde{U}_{14} = 4$ and so on.

Denote by $p_0 = \mathbb{P}(\tilde{U}_{11} = 0, \tilde{U}_{12} = 0, ..., \tilde{U}_{1d} = 0) \geq 0$, $g_{U_{1s}}(z) = \mathbb{E} z^{U_{1s}}$, $g_{\tilde{U}_{1s}}(z) = \mathbb{E} z^{\tilde{U}_{1s}}$ and by $\overrightarrow{U}_i = (U_{i1}, U_{i2}, ..., U_{id}) =$
$$= (\tilde{U}_{i1}, \tilde{U}_{i2}, ..., \tilde{U}_{id}|(\tilde{U}_{i1}, \tilde{U}_{i2}, ..., \tilde{U}_{id}) \not=(0, 0, ..., 0)), i = 1, 2, ...$$

We need to satisfy the no-empty claims requirement in the Cramer-Lundberg model.  Therefore we consider also thinned Poisson process $N = \{N(t): t \geq 0\}$ which counts only the non-empty groups. Due to the independence between $\tilde{N}(t)$ and $\overrightarrow{\tilde{U}}_1, \overrightarrow{\tilde{U}}_2, ..., $ it is well known that $N$ is a HPP process with intensity $\lambda = \tilde{\lambda}[1 - p_0]$ and the random vectors $\overrightarrow{U}_1, \overrightarrow{U}_2, ..., $ are also i.i.d. and independent on $N$. It is easy to obtain the relations between the numerical characteristics of $\overrightarrow{\tilde{U}}_1$ and $\overrightarrow{U}_1$ therefore we will skip the proof of the next result.

{\bf Lemma 4.1} For $s = 1, 2, ..., d$,

$$\mathbb{E} U_{1s} = \frac{\mathbb{E}\tilde{U}_{1s}}{1 - p_0}, \quad g_{U_{1s}}(z) = \frac{g_{\tilde{U}_{1s}}(z)- p_0}{1 - p_0}, \quad \mathbb{D} U_{1s} = \frac{\mathbb{D}\tilde{U}_{1s}}{1 - p_0} - \frac{p_0(\mathbb{E}\tilde{U}_{1s})^2}{(1 - p_0)^2}.$$

\bigskip

Denote by $\overrightarrow{N}_c$ the number of claims of different types that have occurred up to time $t > 0$ without counting the empty groups. Then this is a vector-valued process
$\overrightarrow{N}_c = \{\overrightarrow{N}_c(t): t \geq 0\}$
and its $s$-th coordinate is
$$N_{c,s}(t) = I\{N(t) > 0\} \sum_{i=1}^{N(t)} U_{is}, \,\, s = 1, 2, ..., d.$$
It represents the number of claims of type $s$ that have occurred up to time $t > 0$. It order to be able to count also the empty groups we consider the vector valued process $\overrightarrow{\tilde{N}}_c = \{\overrightarrow{\tilde{N}}_c(t): t \geq 0\}$
$$\tilde{N}_{c,s}(t) = I\{\tilde{N}(t) > 0\} \sum_{i=1}^{\tilde{N}(t)} \tilde{U}_{is}, \,\, s = 1, 2, ..., d.$$
On the realization depicted on Figure \ref{fig:Fig2}, $\tilde{N}_{c,1}(t) \geq 9$, $\tilde{N}_{c,2}(t) \geq 12$, $\tilde{N}_{c,3}(t) \geq 8$ and $\tilde{N}_{c,s}(t) \geq 15$.

Having in mind Lemma 3.1 we obtain that these two processes are stochastically equivalent in the sense of the equality of their finite dimensional distributions (f.d.ds).

{\bf Lemma 4.2} $\tilde{N}_{c} \stackrel{\rm f.d.d.}{=} N_{c}.$

{\it Proof.}  First let us consider the univariate marginals. For all $t \geq 0$ we can check the equality between the corresponding p.g.fs. The p.g.f. of the compound Poisson process is well known. Using its form and Lemma 4.1. we obtain
$$\mathbb{E} z^{I\{N(t) > 0\} \sum_{i=1}^{N(t)} U_{is}} = e^{-\tilde{\lambda}(1 - p_0) t(1 - g_{U_{1s}}(z))} = e^{-\tilde{\lambda}(1 - p_0) t(1 - \frac{g_{\tilde{U}_{1s}}(z)- p_0}{1 - p_0})} =$$
 $$= e^{-\tilde{\lambda}t(1 - g_{\tilde{U}_{1s}}(z))} = \mathbb{E} z^{I\{\tilde{N}(t) > 0\} \sum_{i=1}^{\tilde{N}(t)} \tilde{U}_{is}}, \,\, s = 1, 2, ..., d.$$
Due to the independence and homogeneity of the additive increments of homogeneous Poisson process, the analogous equalities could be proven for the additive increments and consequently for all their finite dimensional distributions.\qed\,.

Particular cases of the counting process  $\overrightarrow{\tilde{N}}_c$ are investigated by many authors. These processes coincide with Compound Poisson processes with non-negative integer valued summands. For $\overrightarrow{\tilde{U}}_1$
\begin{itemize}
  \item Multinomial;
  \item Negative multinomial or in particular multivariate geometric;
  \item Shifted negative binomial and in particular Shifted geometric or
  \item Poisson
\end{itemize}
distributed the time intersections of these processes are investigated in 1962, by G. Smith \cite{SmithPhDThesis}. The particular case, when
the $\overrightarrow{U}_1$ are Shifted multivariate geometric on the set of natural numbers, is investigated in series of papers of Minkova and Balakrishnan \cite{minkova2014type,minkova2014bivariate}.

In order to reduce this model to the Cramer-Lundberg model let us denote the total claim amount within $i$-th group, $i = 1, 2, ...$ by
$$Y_{i} = \sum_{s = 1}^d I\{U_{is} > 0\}\sum_{j = 1}^{U_{is}} \tilde{Y}_{ijs}.$$
These  r.vs. are strictly positive, i.i.d., having common non-lattice c.d.f., with mean
$$\mathbb{E}(Y_1) = \sum_{s = 1}^d \mathbb{E} U_{1s} \mathbb{E} \tilde{Y}_{11s} = \sum_{s = 1}^d \frac{\mathbb{E}\tilde{U}_{1s}}{1 - p_0} \mathbb{E} \tilde{Y}_{11s} < \infty.$$
In case, when the coordinates of $\overrightarrow{\tilde{Y}}_{111}$ are independent
$$\mathbb{D}(Y_1) =  \sum_{s = 1}^d\left[\frac{\mathbb{D}\tilde{U}_{1s}}{1 - p_0} - \frac{p_0(\mathbb{E}\tilde{U}_{1s})^2}{(1 - p_0)^2}\right](\mathbb{E} \tilde{Y}_{11s})^2 + \sum_{s = 1}^d \frac{\mathbb{E} \tilde{U}_{1s}}{1- p_0}\mathbb{D} \tilde{Y}_{11s}\leq \infty.$$
In case when we have only one possible type, i.e. $d = 1$
$$FI\,\,Y_1 =   \left[FI\,\,\tilde{U}_{11} - \frac{p_0\mathbb{E}\tilde{U}_{11}}{1 - p_0}\right]\mathbb{E} \tilde{Y}_{11} + FI\,\, \tilde{Y}_{11}.$$
If we take into account also the empty groups(for example if a claim of a group occur but finally it turns out that there is no real reason to pay for it) it is clear that
$$\tilde{Y}_{i} = \sum_{s = 1}^d I\{\tilde{U}_{is} > 0\}\sum_{j = 1}^{\tilde{U}_{is}} \tilde{Y}_{ijs}, i = 1, 2, ...$$
$\mathbb{P}(\tilde{Y}_{i} = 0) = p_0 \geq 0$ and their mean is
$$\mathbb{E}(\tilde{Y}_1) = \sum_{s = 1}^d \mathbb{E} \tilde{U}_{1s} \mathbb{E} \tilde{Y}_{11s} = \mathbb{E}(Y_1)(1 - p_0) < \infty.$$

Denote the c.d.f. of $\tilde{Y}_{11s}$ with $F_s$. Then for $x \geq 0$,
$$F_{Y_1}(x) = \sum_{u_1 = L_1}^{R_1} ... \sum_{u_d = L_d}^{R_d} \mathbb{P}(\sum_{s = 1}^d I\{u_s > 0\}\sum_{j = 1}^{u_s} \tilde{Y}_{1js} \leq x) \mathbb{P}(U_{11} = u_1, ..., U_{1d} = u_d).$$
In case when the coordinates of $\overrightarrow{\tilde{Y}}_{111}$ are independent we can use the convolution operator and

\begin{equation}\label{cdfY1}
F_{Y_1}(x) = \sum_{u_1 = \tilde{L}_1}^{R_1} ... \sum_{u_d = \tilde{L}_d}^{R_d} F_1^{u_1*}*F_2^{u_2*} *...* F_d^{u_d*}(x)
\end{equation}
$$.\frac{P(\tilde{U}_{11} = u_1, ..., \tilde{U}_{1d} = u_d)}{1 - p_0} - \frac{p_0}{1 - p_0} = \frac{F_{\tilde{Y}_1}(x) - p_0}{1 - p_0},$$
where $0 \leq \tilde{L}_s$ and $\tilde{R}_s$ are correspondingly the lower and the upper end points of the distribution of $\tilde{U}_{1s}$, and $0 \leq L_s$ and $R_s = \tilde{R}_s$ are correspondingly the lower and the upper end points of the distribution of $U_{1s}$, $s = 1, 2, ..., d$.

$\bar{F}_{Y_1}(x) =$
$$=  \sum_{u_1 = \tilde{L}_1}^{R_1} ... \sum_{u_d = \tilde{L}_d}^{R_d} [1-F_1^{u_1*}*F_2^{u_2*} *...* F_d^{u_d*}(x)]\frac{P(\tilde{U}_{11} = u_1, ..., \tilde{U}_{1d} = u_d)}{1 - p_0} =$$
$$= \frac{\bar{F}_{\tilde{Y}_1}(x)}{1 - p_0}.$$
The independence between the coordinates is used here only for making convolution between the c.d.fs of the coordinates. However it is important for the following form of the Laplace - Stieltjes transform
\begin{equation}\label{LStY1GRoup}
l_{Y_1}(s) =  \frac{g_{\tilde{U}_{11}, ..., \tilde{U}_{1d}}(l_{\tilde{Y}_{111}}(s), l_{\tilde{Y}_{112}}(s), ..., l_{\tilde{Y}_{11d}}(s)) - p_0}{1 - p_0} =
\end{equation}
$$= \frac{l_{\tilde{Y}_1}(s) - p_0}{1 - p_0}.$$
And it turns out that the integrated tail distributions of $Y_1$ and $\tilde{Y}_1$ coincide. Their Laplace - Stieltjes transforms are
$$l_{F_{I, Y_1}}(s) = \frac{1}{\mathbb{E}Y_1} \int_0^\infty e^{-sy} \bar{F}_{Y_1}(y)dy  = \frac{1 - l_{Y_1}(s)}{s\mathbb{E}Y_1} = $$
$$= \frac{1 - g_{U_{11}, ..., U_{1d}}(l_{\tilde{Y}_{111}}(s), l_{\tilde{Y}_{112}}(s), ..., l_{\tilde{Y}_{11d}}(s))}{s\sum_{i = 1}^d \mathbb{E}U_{1s}\mathbb{E}\tilde{Y}_{11s}} =$$
$$ = \frac{1 - \frac{g_{\tilde{U}_{11}, ..., \tilde{U}_{1d}}(l_{\tilde{Y}_{111}}(s), l_{\tilde{Y}_{112}}(s), ..., l_{\tilde{Y}_{11d}}(s)) - p_0}{1 - p_0}}{s\sum_{i = 1}^d \mathbb{E}U_{1s}\mathbb{E}\tilde{Y}_{11s}} =$$
$$=  \frac{1 - g_{\tilde{U}_{11}, ..., \tilde{U}_{1d}}(l_{\tilde{Y}_{111}}(s), l_{\tilde{Y}_{112}}(s), ..., l_{\tilde{Y}_{11d}}(s))}{s(1 - p_0)\sum_{i = 1}^d \mathbb{E}U_{1s}\mathbb{E}\tilde{Y}_{11s}} =$$
$$= \frac{1 - g_{\tilde{U}_{11}, ..., \tilde{U}_{1d}}(l_{\tilde{Y}_{111}}(s), l_{\tilde{Y}_{112}}(s), ..., l_{\tilde{Y}_{11d}}(s))}{s\sum_{i = 1}^d \mathbb{E}\tilde{U}_{1s}\mathbb{E}\tilde{Y}_{11s}} = l_{F_{I, \tilde{Y}_1}}(s).$$
The last means that both $Y_1$ and $\tilde{Y}_1$ have equal in distribution equilibrium(or integrated tail) distributions.

The total claim amount up to time $t$ then can be presented by the traditional equality of the Cramer-Lundberg model. It is
$$S(t) = I\{N(t) > 0\} \sum_{i = 1}^{N(t)} Y_{i} = I\{N(t) > 0\} \sum_{i = 1}^{N(t)} \sum_{s = 1}^d I\{U_{is} > 0\}\sum_{j = 1}^{U_{is}}\tilde{Y}_{ijs} $$
$$ \stackrel{\rm d}{=} \sum_{s = 1}^d I\{N_{c,s}(t) > 0\}\sum_{i = 1}^{N_{c,s}(t)}\tilde{Y}_{is} \stackrel{\rm d}{=} \sum_{s = 1}^d I\{\tilde{N}_{c,s}(t) > 0\}\sum_{i = 1}^{\tilde{N}_{c,s}(t)}\tilde{Y}_{is},$$
where $\tilde{Y}_{is}$, $i = 1, 2, ...$ are i.i.d. auxiliary or renumbered random variables in the same probability space, which coincide in distribution with $\tilde{Y}_{11s}$.

Consider the following Risk process
\begin{equation}\label{GeneralModelGroups}
R(t) = u + ct - S(t).
\end{equation}
It is obviously a particular case of the Cramer-Lundberg model. Thus we derive general formulae for the numerical characteristics  and relate  them with probabilities of ruin in this model.

The safety loading is
\begin{equation}\label{rhogroup}
\rho = \frac{c}{\lambda \sum_{s = 1}^d \mathbb{E}U_{1s} \mathbb{E} \tilde{Y}_{1s}} - 1 = \frac{c}{\tilde{\lambda}\sum_{s = 1}^d \mathbb{E}\tilde{U}_{1s} \mathbb{E} \tilde{Y}_{1s}} - 1.
\end{equation}

In general, by (\ref{NPC}) the net profit condition: $\rho > 0$, is equivalent to
\begin{equation}\label{NPCGroups}
c >  \lambda \sum_{s = 1}^d \mathbb{E}U_{1s} \mathbb{E} \tilde{Y}_{1s} \Leftrightarrow c >  \tilde{\lambda} \sum_{s = 1}^d \mathbb{E}\tilde{U}_{1s} \mathbb{E} \tilde{Y}_{1s}
\end{equation}
Equality (\ref{0}) gives immediately
\begin{equation}\label{0Groups}
\psi(0) =  \frac{\lambda}{c}\sum_{s = 1}^d \mathbb{E}U_{1s} \mathbb{E} \tilde{Y}_{1s}, \,\, \delta(0) = 1 - \frac{\lambda}{c}\sum_{s = 1}^d \mathbb{E}U_{1s} \mathbb{E} \tilde{Y}_{1s}.
\end{equation}

In case when the coordinates of $\overrightarrow{\tilde{Y}}_{111}$ are independent, the distribution  of the deficit at the time of ruin, given that ruin with initial capital zero has the following distribution
$$\mathbb{P}(-R(\tau(0)+)\leq x| \tau(0) < \infty) = F_{I, Y_1}(x) = \frac{1}{\mathbb{E}Y_1}\int_0^x \bar{F}_{Y_1}(y)dy  = $$
$$= \frac{1}{\mathbb{E}Y_1} \sum_{u_1 = L_1}^{R_1} ... \sum_{u_d = L_d}^{R_d}P(U_{11} = u_1, ..., U_{1d} = u_d)\int_0^x  [1- F_1^{u_1*} *...* F_d^{u_d*}(y)]dy = $$
$$= \frac{1}{(1 - p_0)\mathbb{E}Y_1} \sum_{u_1 = \tilde{L}_1}^{R_1} ... \sum_{u_d = \tilde{L}_d}^{R_d} P(\tilde{U}_{11} = u_1, ..., \tilde{U}_{1d} = u_d)\int_0^x [1- F_1^{u_1*}*... $$
\hfill $* F_d^{u_d*}(y)]dy =$
\begin{equation}\label{intyegratedtailY1Groups}
= \frac{1}{\mathbb{E}\tilde{Y}_1} \int_0^x [1-F_{\tilde{Y}_1}(y)]dy = F_{I, \tilde{Y}_1}(x).
\end{equation}
The independence between the coordinates is used here only for making convolution between the c.d.fs of the coordinates.

We receive for its expectation

$\mathbb{E}(-R(\tau(0)+)| \tau(0) < \infty) = $
$$  = \frac{\sum_{s = 1}^d \mathbb{E}\tilde{U}_{1s}\mathbb{E} \tilde{Y}_{11s}}{2(1 - p_0)} +
 \frac{\sum_{s = 1}^d\left[\mathbb{D}\tilde{U}_{1s} - \frac{p_0}{1 - p_0}(\mathbb{E}\tilde{U}_{1s})^2\right](\mathbb{E} \tilde{Y}_{11s})^2 + \sum_{s = 1}^d \mathbb{E} \tilde{U}_{1s}\mathbb{D} \tilde{Y}_{11s}}{2\left(\sum_{s = 1}^d \mathbb{E}\tilde{U}_{1s}\mathbb{E} \tilde{Y}_{11s} \right)}.$$
 For $d = 1$ one can easily obtain
 $$ \mathbb{E}(-R(\tau(0)+)| \tau(0) < \infty) = \frac{1}{2} \left[(FI \tilde{U}_{11} + \mathbb{E}{\tilde{U}_{11}})\mathbb{E} \tilde{Y}_{11} + FI \tilde{Y}_{11}\right].$$

For $\mathbb{D}\tilde{Y}_{1s} < \infty$ and $\mathbb{D} U_{1s} < \infty$, $d = 1, 2, ..., d$ the equation (\ref{tau0}) gives us immediately

$\mathbb{E}(\tau(0)|\tau(0) < \infty) =$
$$=  \frac{\sum_{s = 1}^d\left[\mathbb{D}\tilde{U}_{1s} - \frac{p_0}{1 - p_0}(\mathbb{E}\tilde{U}_{1s})^2\right](\mathbb{E} \tilde{Y}_{11s})^2 + \sum_{s = 1}^d \mathbb{E} \tilde{U}_{1s}\mathbb{D} \tilde{Y}_{11s}+ \frac{1}{1 - p_0}\left(\sum_{s = 1}^d \mathbb{E}\tilde{U}_{1s}\mathbb{E} \tilde{Y}_{11s} \right)^2}{2\left(c - \tilde{\lambda}\sum_{s = 1}^d \mathbb{E}\tilde{U}_{1s} \mathbb{E} \tilde{Y}_{1s}\right)\left(\sum_{s = 1}^d \mathbb{E}\tilde{U}_{1s}\mathbb{E} \tilde{Y}_{11s} \right)}.$$
And for $d = 1$

$\mathbb{E}(\tau(0)|\tau(0) < \infty) =$
 $$  =\frac{(FI U_{11} + \mathbb{E}U_{11})\mathbb{E} \tilde{Y}_{11s} + FI \tilde{Y}_{11}}{2\left(c - \lambda\mathbb{E}U_{11} \mathbb{E} \tilde{Y}_{11}\right)} = \frac{(FI\tilde{U}_{11} + \mathbb{E}\tilde{U}_{11})\mathbb{E} \tilde{Y}_{11} + FI \tilde{Y}_{11}}{2\left(c - \tilde{\lambda}\mathbb{E}\tilde{U}_{11} \mathbb{E} \tilde{Y}_{11}\right)} .
 $$

Formulas (\ref{jointdeficite}) and (\ref{intyegratedtailY1Groups}) imply that the joint distribution of the severity of (deficit at) ruin and the risk surplus just before the ruin with initial capital zero is:

$\mathbb{P}(-R(\tau(0)+) > x, R(\tau(0)-) > y|\tau(0) < \infty) =$
$$ = \frac{1}{\sum_{s = 1}^d \mathbb{E}\tilde{U}_{1s}\mathbb{E} \tilde{Y}_{11s}} \sum_{u_1 = \tilde{L}_1}^{R_1} ... \sum_{u_d = \tilde{L}_d}^{R_d} P(\tilde{U}_{11} = u_1, ..., \tilde{U}_{1d} = u_d) $$
$$\int_{x+y}^\infty [1-F_1^{u_1*}*F_2^{u_2*} *...* F_d^{u_d*}(z)]dz.$$

The distribution of the claim causing ruin (\ref{theclaimcausingruin}) in this case is

$\mathbb{P}(R(\tau(0)-)-R(\tau(0)+) \leq x| \tau(0) < \infty) =$
$$= \frac{1}{\sum_{s = 1}^d \mathbb{E}\tilde{U}_{1s}\mathbb{E} \tilde{Y}_{11s}} \sum_{u_1 = \tilde{L}_1}^{R_1} ... \sum_{u_d = \tilde{L}_d}^{R_d}P(\tilde{U}_{11} = u_1, ..., \tilde{U}_{1d} = u_d) \int_0^x y d F_1^{u_1*}*$$
\hfill $*F_2^{u_2*} *...* F_d^{u_d*}(x).$

For $u > 0$, using formula (\ref{deltaPrim}), (\ref{0}) and the form of the c.d.f. (\ref{cdfY1}) we obtain

$$\delta'(u) = \frac{\lambda}{c}\left[\delta(u) - \sum_{u_1 = \tilde{L}_1}^{R_1} ... \sum_{u_d = \tilde{L}_d}^{R_d} \frac{P(\tilde{U}_{11} = u_1, ..., \tilde{U}_{1d} = u_d)}{1 - p_0}\right.$$
\begin{equation}\label{deltaPrimGroups}
\left. \int_0^u \delta(u - y) d  F_1^{u_1*} *...* F_d^{u_d*}(x)\right]
\end{equation}

Equations (\ref{renewaldelta}) and (\ref{renewalpsi}) imply that the solution of (\ref{deltaPrimGroups}) satisfy the defective renewal equation
$$ \delta(u) = 1 - \frac{\lambda}{c}\sum_{s = 1}^d \mathbb{E}U_{1s} \mathbb{E} \tilde{Y}_{1s} +  \frac{\lambda}{c}\sum_{u_1 = L_1}^{R_1} ... \sum_{u_d = L_d}^{R_d}P(U_{11} = u_1, ..., U_{1d} = u_d)$$
$$\int_0^u \delta(u-y) [1- F_1^{u_1*} *...* F_d^{u_d*}(y)]dy.$$
and
$$ \psi(u) = \frac{\tilde{\lambda}}{c}\sum_{u_1 = \tilde{L}_1}^{R_1} ... \sum_{u_d = \tilde{L}_d}^{R_d}P(\tilde{U}_{11} = u_1, ..., \tilde{U}_{1d} = u_d)\int_0^x  [1- F_1^{u_1*} *...* F_d^{u_d*}(y)]dy - $$
$$- \frac{\tilde{\lambda}}{c}\sum_{u_1 = \tilde{L}_1}^{R_1} ... \sum_{u_d = \tilde{L}_d}^{R_d}P(\tilde{U}_{11} = u_1, ..., \tilde{U}_{1d} = u_d)\int_0^u \psi(u-y) [1- F_1^{u_1*} *...* F_d^{u_d*}(y)]dy.$$

%{

%Then the formulae (\ref{GeneralPsi}) and (\ref{deltageneral}) imply
%\begin{equation}\label{PloachekGroupsPsi}
%\psi(u) =   \left(1 - \psi(0)\right) \sum_{i = 1}^\infty \psi^i(0) [1 - F_I^{i*}(u)], \,\, u \geq 0.
%\end{equation}
%\begin{equation}\label{PloachekGroupsdelta}
%\delta(u)  = \left(1 - \psi(0)\right) \sum_{i = 0}^\infty \psi^i(0) F_I^{i*}(u),  \,\, u \geq 0.
%\end{equation}

%}

In case when the coordinates for $\overrightarrow{\tilde{Y}}_{1}$ are independent the equations (\ref{LStY1GRoup}), (\ref{LStTrdeltaGeneral}), and (\ref{rhogroup}) allow us to obtain a general formula for the Laplace - Stieltjes transform of $\delta(u)$
$$l_{\delta}(s) =  \frac{1 - \frac{\lambda}{c}\sum_{s = 1}^d \mathbb{E}U_{1s} \mathbb{E} \tilde{Y}_{1s}}{1 - \frac{\lambda}{cs} \left(1 - \frac{g_{\tilde{U}_{11}, ..., \tilde{U}_{1d}}(l_{\tilde{Y}_{111}}(s), l_{\tilde{Y}_{112}}(s), ..., l_{\tilde{Y}_{11d}}(s)) - p_0}{1 - p_0} )\right)} =$$
$$= \frac{1 - \frac{\lambda}{c}\sum_{s = 1}^d \mathbb{E}U_{1s} \mathbb{E} \tilde{Y}_{1s}}{1 - \frac{\tilde{\lambda}}{cs} \left[1 - g_{\tilde{U}_{11}, ..., \tilde{U}_{1d}}(l_{\tilde{Y}_{111}}(s), l_{\tilde{Y}_{112}}(s), ..., l_{\tilde{Y}_{11d}}(s)\right]}.$$

By formula (\ref{G0y})

$G(0, y) =$
$$ = \frac{\tilde{\lambda}}{c}\sum_{u_1 = \tilde{L}_1}^{R_1} ... \sum_{u_d = \tilde{L}_d}^{R_d} P(\tilde{U}_{11} = u_1, ..., \tilde{U}_{1d} = u_d) \int_0^x [1-F_1^{u_1*}*...* F_d^{u_d*}(y)]dy$$

Moreover
\begin{equation*}\label{GGroups}
G(u, y) = \frac{\tilde{\lambda}}{c}\left[\sum_{u_1 = \tilde{L}_1}^{R_1} ... \sum_{u_d = \tilde{L}_d}^{R_d}P(\tilde{U}_{11} = u_1, ..., \tilde{U}_{1d} = u_d) \int_0^u G(u - x, y) [1-\right.$$
$$ - F_1^{u_1*}*...* F_d^{u_d*}(x)]dx + \sum_{u_1 = \tilde{L}_1}^{R_1} ... \sum_{u_d = \tilde{L}_d}^{R_d} P(\tilde{U}_{11} = u_1,..., \tilde{U}_{1d} = u_d) $$
$$ \left.  \int_u^{u + y}[1-F_1^{u_1*}*F_2^{u_2*} *...* F_d^{u_d*}(x)]dx\right], \,\, u \geq 0.
\end{equation*}

%Consider again the case, when the coordinates of $\overrightarrow{\tilde{Y}}_{11}$ are independent. For $r > 0$, such that $g_{\tilde{U}_{11}, ..., \tilde{U}_{1d}}(l_{\tilde{Y}_{111}}(-r), l_{\tilde{Y}_{112}}(-r), ..., l_{\tilde{Y}_{11d}}(-r)) < \infty$, denote by
%$$g_{gr}(r) = \tilde{\lambda}[g_{\tilde{U}_{11}, ..., \tilde{U}_{1d}}(l_{\tilde{Y}_{111}}(-r), l_{\tilde{Y}_{112}}(-r), ..., l_{\tilde{Y}_{11d}}(-r)) - 1] - cr.$$ Then, the process
%$$\{e^{-r[R(t)-u]-t\left[\tilde{\lambda}[g_{\tilde{U}_{11}, ..., \tilde{U}_{1d}}(l_{\tilde{Y}_{111}}(-r), l_{\tilde{Y}_{112}}(-r), ..., l_{\tilde{Y}_{11d}}(-r)) - 1] - cr\right]}: t \geq 0\}$$
%is a martingale and for all $T > 0$ and $r > 0$,
%$$
%\nonumber \psi(u, T) \leq  \frac{e^{-ru}}{\mathbb{E}[e^{-g_{gr}(r)\tau(u)}| \tau(u) < T]} \leq e^{-ru}\sup_{y \in [0, T]}e^{yg_{gr}(r)}.
%$$

For this model, the small claim condition  means that there exists the Cramer-Lundberg exponent $\epsilon > 0$:
$$g_{\tilde{U}_{11}, ..., \tilde{U}_{1d}}(l_{\tilde{Y}_{111}}(-\epsilon), l_{\tilde{Y}_{112}}(-\epsilon), ..., l_{\tilde{Y}_{11d}}(-\epsilon)) = \frac{c}{\tilde{\lambda}}\epsilon + 1.$$
and for this $\epsilon$ the inequality (\ref{psiLundbergexponent}), that is $\psi(u) \leq e^{-\epsilon u}$ is satisfied and choosing $\alpha \in (0, 1)$ we can find appropriate premium income rate $c$, such that $\psi(u) \leq \alpha$.

 If additionally
$$ \sum_{u_1 = \tilde{L}_1}^{R_1} ... \sum_{u_d = \tilde{L}_d}^{R_d} \mathbb{P}(\tilde{U}_{11} = u_1, ..., \tilde{U}_{1d} = u_d) \int_0^\infty x e^{\epsilon x} [1-F_1^{u_1*}*...* F_d^{u_d*}(x)] dx < \infty,$$
then the following Cramer-Lundberg approximation of the probability of ruin holds
$$\lim_{u \to \infty} e^{\epsilon u}\psi(u) = \frac{c - \tilde{\lambda}\sum_{s = 1}^d \mathbb{E}U_{1s} \mathbb{E} \tilde{Y}_{1s}}{\epsilon\tilde{\lambda}}\sum_{u_1 = \tilde{L}_1}^{R_1} ... \sum_{u_d = \tilde{L}_d}^{R_d} \mathbb{P}(\tilde{U}_{11} = u_1, ...,  \tilde{U}_{1d} = u_d)$$
$$\int_0^\infty x e^{\epsilon x} [1-F_1^{u_1*}*...* F_d^{u_d*}(x)] dx.$$

{If $Y_1$ belongs to the class $D$ of dominatedly varying distributions then, given the net profit condition is satisfied, we can use the asymptotic (\ref{AsymptoticHeavy}) and obtain that
$$
\lim_{u \to \infty} \frac{\psi(u)}{\sum_{u_1 = \tilde{L}_1}^{R_1} ... \sum_{u_d = \tilde{L}_d}^{R_d} P(\tilde{U}_{11} = u_1, ..., \tilde{U}_{1d} = u_d) \int_{u}^\infty [1-F_1^{u_1*}*...* F_d^{u_d*}(z)]dz} =$$
$$= \frac{\tilde{\lambda}}{c - \tilde{\lambda}\sum_{s = 1}^d \mathbb{E}\tilde{U}_{1s} \mathbb{E} \tilde{Y}_{1s}}. $$
To this class belong the most of the claim size distributions that can be met in practice, e.g. Pareto distribution, log-gamma, log-normal, Burr, Weibull with parameter less than 1, Banktander-type-I and type II. For more extensive discussion of the distributions from class $D$ see Embrechts, Klueppelberg, Mikosch \cite{embrechts2013modelling}.}

\section{The case when the counting group processes are independent homogeneous Poisson processes}
In this section we suppose that the numbers of groups of claims that have occurred up to time $t > 0$, inclusively possibly empty groups are independent homogeneous Poisson processes and we reduce this  model to the model discussed in Section 4.

Denote by $\tilde{N}_1(t), \tilde{N}_2(t), ..., \tilde{N}_d(t)$ the numbers of groups of different types $1, 2, ..., d$ that have occurred up to time $t > 0$. See Figure \ref{fig:Fig3} for $d = 4$. The first line depicts the arrival process of the groups of claims of the first type. The second line is for the moments of arrivals of the groups of the second type. In analogous way for the third and the fourth line. The last line represents the merging of these arrival processes.

\begin{center}
\includegraphics[width=.9\textwidth]{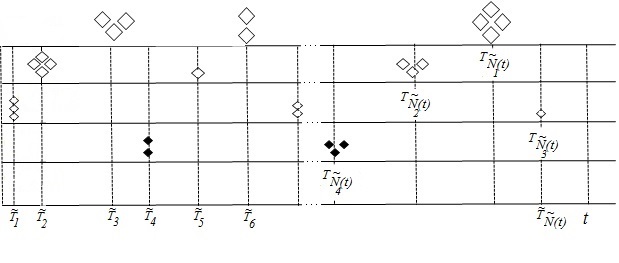}

{\small Fig. 3. The claim arrival processes, when the counting group processes are independent. \label{fig:Fig3}}
\end{center}

 Assume that these processes are independent and in order to obtain a particular case of the Cramer-Lundberg model we assume that $\tilde{N}_s = \{\tilde{N}_s(t): t \geq 0\}$ are independent HPPs with intensity $\tilde{\lambda}_i$, $s = 1, 2, ..., d$. Because of these processes are stochastically continuous and independent simultaneous jumps are almost sure not possible. Therefore in this model it is almost sure not possible groups of two or more types of claims to arrive simultaneously.   From the theory of Poisson processes (see e.g. \cite{DaleyVereJones}) it is well known that the merging of $\tilde{N}_1(t), \tilde{N}_2(t), ..., \tilde{N}_d(t)$ is a $\tilde{N} \sim HPP(\tilde{\lambda}_1 + \tilde{\lambda}_2 + ... +\tilde{\lambda}_d)$. Obviously this and the previous model coincide for $d = 1$.

  Let $\overrightarrow{\tilde{U}}_i$ be the number of claims within $i$-th group. The random vectors
$\overrightarrow{\tilde{U}}_i = (\tilde{U}_{i1}, \tilde{U}_{i2}, ..., \tilde{U}_{id})$, $i = 1, 2, ...$
 are i.i.d.  Their coordinates are non-negative integer valued r.vs. The $s$-th coordinate represents the number of claims within $i$-th group and of type $s = 1, 2, ..., d$. On  Figure \ref{fig:Fig3} $d = 4$ and for the outcome of the experiment depicted there the vector $\overrightarrow{\tilde{U}}_1$ has coordinates $\tilde{U}_{11} = 0$, $\tilde{U}_{12} = 0$, $\tilde{U}_{13} = 3$, $\tilde{U}_{14} = 0$. The vector $\overrightarrow{\tilde{U}}_2$ has coordinates $\tilde{U}_{11} = 0$, $\tilde{U}_{12} = 4$, $\tilde{U}_{13} = 0$, $\tilde{U}_{14} = 0$ and so on.

 We would like to show that this model is a particular case of the model discussed in the previous sections and therefore it is a particular case of the Cramer-Lundberg model. Analogous reduction to the following is done e.g. in {\cite{A1, A2, A3}} among others. In order not to lose information about the initial processes and having the merging of these Poisson processes to be able to describe the type of the group which arrives it time $ t> 0$ and is counted by $\tilde{N}$, for any fixed $i = 1, 2, ...$, we denote by $A_{i1}, A_{i2}, ..., A_{id}$ stochastically equivalent and independent partitions of the sample space $\Omega$. More precisely these event are mutually exclusive and for all $i = 1, 2, ...$
 $$\bigcup_{s = 1}^d A_{is} = \Omega, \quad \mathbb{P}(A_{is}) = \frac{\tilde{\lambda}_s}{\tilde{\lambda}_1 + \tilde{\lambda}_2 + ... + \tilde{\lambda}_d}, \,\, s = 1, 2, ..., d.$$

Denote by
$$p_0 = \mathbb{P}(I\{A_{11}\}\tilde{U}_{i1} = 0, I\{A_{12}\}\tilde{U}_{i2} = 0, ..., I\{A_{1d}\}\tilde{U}_{id} = 0) = $$
$$= \frac{\sum_{s = 1}^d \mathbb{P}(\tilde{U}_{1s} = 0)\tilde{\lambda}_s}{\tilde{\lambda}_1 + \tilde{\lambda}_2 + ... + \tilde{\lambda}_d}$$
the probability that given a group arrived it to be empty. In order to satisfy the non-empty claim requirement in the Cramer-Lundberg model we consider also the random vectors

$\overrightarrow{U}_i = (U_{i1}, U_{i2}, ..., U_{id}) = (I\{A_{i1}\}\tilde{U}_{i1}, I\{A_{i2}\}\tilde{U}_{i2}, ..., I\{A_{id}\}\tilde{U}_{id}|$

  \hfill $|(I\{A_{i1}\}\tilde{U}_{i1}, I\{A_{i2}\}\tilde{U}_{i2}, ..., I\{A_{id}\}\tilde{U}_{id}) \not=(0, 0, ..., 0))$

\noindent for $i = 1, 2, ...$ which count the claims in the simultaneously non-empty groups.
Then the random vectors $\overrightarrow{U}_1, \overrightarrow{U}_2, ..., $ are  i.i.d. Their distribution and numerical characteristics are described in the following lemma.

\medskip

{\bf Lemma 5.1} Let $0 < L_s$ and $R_s$ be correspondingly the lower and the upper end points of the distribution of $U_{1s}$, $s = 1, 2, ..., d$. For $s = 1, 2, ..., d$ and $u_s \in \mathbb{N}$,
$$\mathbb{P}(U_{11} = 0, ..., U_{1s-1} = 0, U_{1s} = u_s, U_{1s+1} = 0, ..., U_{1d} = 0) =$$
$$= \frac{\tilde{\lambda}_s\mathbb{P}(\tilde{U}_{1s} = u_s)}{(\tilde{\lambda}_1 + ... + \tilde{\lambda}_d)(1-p_0)}.$$
  And the last probability is equal to zero in other cases.
$$\mathbb{E} U_{1s} = \frac{\tilde{\lambda_s}\mathbb{E}\tilde{U}_{1s}}{\sum_{j = 1}^d \tilde{\lambda}_j[1-\mathbb{P}(\tilde{U}_{1j} = 0)]}, \quad g_{U_{1s}}(z) = 1 - \frac{\tilde{\lambda}_s[1 - g_{\tilde{U}_{1s}}(z)]}{\tilde{\lambda}_1 + \tilde{\lambda}_2 + ... + \tilde{\lambda}_d},$$
 $$\mathbb{D} U_{1s} = \frac{\tilde{\lambda}_s\{\mathbb{D}\tilde{U}_{1s}+(\mathbb{E}\tilde{U}_{1s})^2\}}{(1 - p_0)(\tilde{\lambda}_1 + \tilde{\lambda}_2 + ... + \tilde{\lambda}_d)} - \frac{\tilde{\lambda}_s^2(\mathbb{E}\tilde{U}_{1s})^2}{(1 - p_0)^2(\tilde{\lambda}_1 + \tilde{\lambda}_2 + ... + \tilde{\lambda}_d)^2}.$$
$$g_{U_{11}, ..., U_{1d}}(z_1, ..., z_d) =  \frac{\sum_{s=1}^d \tilde{\lambda}_sg_{\tilde{U}_{1s}}(z)}{(1 - p_0)(\tilde{\lambda}_1 + \tilde{\lambda}_2 + ... + \tilde{\lambda}_d)} - \frac{p_0}{(1 - p_0)}.$$

{\bf Proof:} $\mathbb{P}(U_{11} = 0, ..., U_{1s-1} = 0, U_{1s} = u_s, U_{1s+1} = 0, ..., U_{1d} = 0) =$
 $$= \mathbb{P}[I\{A_{11}\}\tilde{U}_{11} = 0,...,I\{A_{1s-1}\}\tilde{U}_{1s-1} = 0,I\{A_{1s}\}\tilde{U}_{1s} = u_s,$$
 $$I\{A_{1s+1}\}\tilde{U}_{1s+1} = 0,..., I\{A_{1d}\}\tilde{U}_{id}= 0|(I\{A_{11}\}\tilde{U}_{11}, ...,$$
  \hfill $I\{A_{1d}\}\tilde{U}_{1d}) \not=(0, 0, ..., 0)] =$
  $$= \frac{\mathbb{P}[I\{A_{11}\}\tilde{U}_{11} = 0,...,I\{A_{1s-1}\}\tilde{U}_{1s-1} = 0,I\{A_{1s}\}\tilde{U}_{1s} = u_s,}{\mathbb{P}[(I\{A_{11}\}\tilde{U}_{11}, I\{A_{12}\}\tilde{U}_{12}, ..., I\{A_{1d}\}\tilde{U}_{1d}) \not=(0, 0, ..., 0)]}$$
  $$\frac{I\{A_{1s+1}\}\tilde{U}_{1s+1} = 0,..., I\{A_{1d}\}\tilde{U}_{id}= 0]}{}=$$
  $$= \frac{\mathbb{P}(A_{s1})\mathbb{P}(\tilde{U}_{1s} = u_s)}{1-p_0} = \frac{\tilde{\lambda}_s\mathbb{P}(\tilde{U}_{1s} = u_s)}{(\tilde{\lambda}_1 + ... + \tilde{\lambda}_d)(1-p_0)}.$$
  And the last probability is equal to zero in other cases.
$$\mathbb{E} U_{1s} = \frac{\tilde{\lambda_s}\mathbb{E}\tilde{U}_{1s}}{(\tilde{\lambda}_1 + \tilde{\lambda}_2 + ... + \tilde{\lambda}_d)(1 - p_0)} = \frac{\tilde{\lambda_s}\mathbb{E}\tilde{U}_{1s}}{\sum_{j = 1}^d \tilde{\lambda}_j[1-\mathbb{P}(\tilde{U}_{1j} = 0)]},$$
 $$g_{U_{1s}}(z) = \frac{g_{I\{A_{1s}\}\tilde{U}_{1s}}(z)- p_0}{1 - p_0} = 1 - \frac{\tilde{\lambda}_s[1 - g_{\tilde{U}_{1s}}(z)]}{\tilde{\lambda}_1 + \tilde{\lambda}_2 + ... + \tilde{\lambda}_d},$$
 $$\mathbb{D} U_{1s} = \frac{\mathbb{D}[I\{A_{1s}\}\tilde{U}_{1s}]}{1 - p_0} - \frac{p_0(\mathbb{E}[I\{A_{1s}\}\tilde{U}_{1s}])^2}{(1 - p_0)^2} =$$
%$$=\frac{\tilde{\lambda}_s\mathbb{E}[\tilde{U}_{1s}]^2}{(1 - p_0)(\tilde{\lambda}_1 + \tilde{\lambda}_2 + ... + \tilde{\lambda}_d)} - \frac{\tilde{\lambda}_s^2(\mathbb{E}\tilde{U}_{1s})^2}{(1 - p_0)^2(\tilde{\lambda}_1 + \tilde{\lambda}_2 + ... + \tilde{\lambda}_d)^2}=$$
$$=\frac{\tilde{\lambda}_s\{\mathbb{D}\tilde{U}_{1s}+(\mathbb{E}\tilde{U}_{1s})^2\}}{(1 - p_0)(\tilde{\lambda}_1 + \tilde{\lambda}_2 + ... + \tilde{\lambda}_d)} - \frac{\tilde{\lambda}_s^2(\mathbb{E}\tilde{U}_{1s})^2}{(1 - p_0)^2(\tilde{\lambda}_1 + \tilde{\lambda}_2 + ... + \tilde{\lambda}_d)^2}.$$
$$g_{U_{11}, ..., U_{1d}}(z_1, ..., z_d) =  \frac{\sum_{s=1}^d\tilde{\lambda}_s[g_{\tilde{U}_{1s}}(z)-p_0]}{(1 - p_0)(\tilde{\lambda}_1 + \tilde{\lambda}_2 + ... + \tilde{\lambda}_d)} =$$
$$= \frac{\sum_{s=1}^d \frac{\tilde{\lambda}_s}{\tilde{\lambda}_1 + \tilde{\lambda}_2 + ... + \tilde{\lambda}_d}g_{\tilde{U}_{1s}}(z)-p_0}{1 - p_0}=\frac{\sum_{s=1}^d \tilde{\lambda}_sg_{\tilde{U}_{1s}}(z)}{(1 - p_0)(\tilde{\lambda}_1 + \tilde{\lambda}_2 + ... + \tilde{\lambda}_d)} - \frac{p_0}{(1 - p_0)}.$$
\hfill \qed\,.

Assume that the claim counting process $\{\tilde{N}_c(t): t \geq 0\}$, inclusively empty claims, is such that
 $$ (\tilde{N}_{c,1}(t),..., \tilde{N}_{c,d}(t)) = \left(I\{\tilde{N}_{1}(t) > 0\}\sum_{i = 1}^{\tilde{N}_{1}(t)} \tilde{U}_{i1}, ..., I\{\tilde{N}_{d}(t) > 0\}\sum_{i = 1}^{\tilde{N}_{d}(t)} \tilde{U}_{id}\right). $$

Assume independence between $\tilde{N}(t)$ and $\overrightarrow{\tilde{U}}_1, \overrightarrow{\tilde{U}}_2, ...$.

In order to satisfy non-empty requirement of the claims $Y_1, Y_2, ...$ in the Cramer-Lundberg model we consider also the Poisson processes $N = \{N(t): t \geq 0\}$ which counts the merging of $\tilde{N}_1(t), \tilde{N}_2(t), ..., \tilde{N}_d(t)$ and does not count empty groups. The resulting process $N$ is thinned process of $\tilde{N}$. More precisely $N \sim HPP[(\tilde{\lambda}_1 + \tilde{\lambda}_2 + ... + \tilde{\lambda}_d)(1 - p_0)].$
 It is independent on the random vectors $\overrightarrow{U}_1, \overrightarrow{U}_2, ..$. Then the number of non-empty claims up to time $t > 0$ is described by the vector-valued processes
$\overrightarrow{N}_c = \{\overrightarrow{N}_c(t): t \geq 0\}$
$$ (N_{c,1}(t),..., N_{c,d}(t)) = \left(I\{N(t) > 0\}\sum_{i = 1}^{N(t)} U_{i1}, ..., I\{N(t) > 0\}\sum_{i = 1}^{N(t)} U_{id}\right). $$
Analogously to Lemma 4.2. we obtain that the vector valued processes $\tilde{N}_c$ and $N_c$ coincide in the sense of the f.d.ds.

{\bf Lemma 5.2.} $(N_{c,1},..., N_{c,d})   \stackrel{\rm f.d.d.}{=}  (\tilde{N}_{c,1},..., \tilde{N}_{c,d}) $.

{\it Proof.} Let $t > 0$.

$\mathbb{E} \left[z_1^{\tilde{N}_{c,1}(t)}...z_d^{\tilde{N}_{c,d}(t)}\right] = $
\begin{eqnarray*}
    &=& \mathbb{E} \left[z_1^{I\{\tilde{N}_{1}(t) > 0\}\sum_{i = 1}^{\tilde{N}_{1}(t)} \tilde{U}_{i1}}...z_d^{I\{\tilde{N}_{d}(t) > 0\}\sum_{i = 1}^{\tilde{N}_{d}(t)} \tilde{U}_{id}}\right]  \\
    &=& e^{-\tilde{\lambda}_1 t[1 - g_{\tilde{U}_{11}}(z_1)]}...e^{-\tilde{\lambda}_d t[1 - g_{\tilde{U}_{1d}}(z_d)]} \\
    &=& e^{-(\tilde{\lambda}_1 + ... + \tilde{\lambda}_d)(1 - p_0)t\left\{\frac{1 - \frac{\tilde{\lambda}_1}{\tilde{\lambda}_1 + \tilde{\lambda}_2 + ... + \tilde{\lambda}_d}g_{\tilde{U}_{11}}(z_1) + ... + \frac{\tilde{\lambda}_d}{\tilde{\lambda}_1 + \tilde{\lambda}_2 + ... + \tilde{\lambda}_d}g_{\tilde{U}_{1d}}(z_d) \pm p_0}{1 - p_0}\right\}} \\
   &=& e^{-(\tilde{\lambda}_1 + ... + \tilde{\lambda}_d)(1 - p_0)t\left\{1 -\frac{\frac{\tilde{\lambda}_1}{\tilde{\lambda}_1 + \tilde{\lambda}_2 + ... + \tilde{\lambda}_d}g_{\tilde{U}_{11}}(z_1) + ... + \frac{\tilde{\lambda}_d}{\tilde{\lambda}_1 + \tilde{\lambda}_2 + ... + \tilde{\lambda}_d}g_{\tilde{U}_{1d}}(z_d)- p_0}{1 - p_0}\right\}} \\
   &=& l_{N(t)}\left\{\frac{\sum_{s=1}^d \frac{\tilde{\lambda}_s}{\tilde{\lambda}_1 + \tilde{\lambda}_2 + ... + \tilde{\lambda}_d}g_{\tilde{U}_{1s}}(z_s)-p_0}{1 - p_0}\right\} \\
   &=& g_{N(t)} (g_{U_{11},...,U_{1d}}(z_1, ..., z_d))\\
    &=& g_{I_{N(t) > 0}\sum_{i = 1}^{N(t)} U_{i1},...,I_{N(t) > 0}\sum_{i = 1}^{N(t)} U_{id}}(z_1, ...,z_d) = \mathbb{E} \left[z_1^{N_{c,1}(t)}...z_d^{N_{c,d}(t)}\right] .
\end{eqnarray*}
Due to the independence and homogeneity of the additive increments of homogeneous Poisson process, the analogous equalities could be proven for the additive increments and consequently for all finite dimensional distributions.\qed\,.

Denote by $\overrightarrow{\tilde{Y}}_{ij}$ the $j$-th claim size within $i$-th group, $i, j = 1, 2, ...$. Due to our assumption that
$d$ types of claims are possible, $\overrightarrow{\tilde{Y}}_{ij}$ is a random vector,
$\overrightarrow{\tilde{Y}}_{ij} = (\tilde{Y}_{ij1}, \tilde{Y}_{ij2}, ...., \tilde{Y}_{ijd}),$
where $\tilde{Y}_{ijs}$ represents $j$-th claim size of type $s$, within $i$-th group, $s = 1, 2, ..., d$. By assumption $F_{\overrightarrow{\tilde{Y}}_{ij}} (0, 0, ..., 0) = 0$. The claim amounts within a group are independent on the inter-arrival times between groups. The claim sizes of a given type are assumed to be i.i.d. Therefore $\tilde{Y}_{ijs} \stackrel{\rm d}{=} \tilde{Y}_{11s}$ for all $i, j \in \mathbb{N}$ and in cases when it is clear we will skip the first index.

The total claim amount up to time $t$ is
$$S(t) = I\{N(t) > 0\} \sum_{i = 1}^{N(t)} Y_{i} = I\{N(t) > 0\} \sum_{i = 1}^{N(t)} \sum_{s = 1}^d I\{U_{is} > 0\}\sum_{j = 1}^{U_{is}}\tilde{Y}_{ijs} $$
$$ \stackrel{\rm d}{=} \sum_{s = 1}^d I\{N_{c,s}(t) > 0\}\sum_{i = 1}^{N_{c,s}(t)}\tilde{Y}_{is} \stackrel{\rm d}{=} \sum_{s = 1}^d I\{\tilde{N}_{c,s}(t) > 0\}\sum_{i = 1}^{\tilde{N}_{c,s}(t)}\tilde{Y}_{is},$$
where $\tilde{Y}_{is}$, $i = 1, 2, ...$ are i.i.d. auxiliary or renumbered random variables in the same probability space, which coincide in distribution with $\tilde{Y}_{11s}$.

Now we reduced the considered process to the Cramer-Lundberg Risk model
\begin{equation}\label{GeneralModelGroupsIndependent}
R(t) = u + ct - S(t).
\end{equation}
where $c > 0$ is the premium income rate and $u \geq 0$ is the initial capital.

Let us now, using the results for the Cramer-Lundberg model to derive general formulae for the numerical characteristics  and related with them probabilities for ruin in this model.

{\bf Theorem~5.1}. The risk process defined in (\ref{GeneralModelGroupsIndependent}) is a particular case of Cramer-Lundberg risk process with initial capital $u$ and premium income rate $c$.
Moreover
\begin{description}
  \item[a.)] the intensity of the corresponding homogeneous Poisson counting process described in (C3) and (C6) is
  $$\lambda = (\tilde{\lambda}_1 + \tilde{\lambda}_2 + ... + \tilde{\lambda}_d)(1 - p_0) = \sum_{s = 1}^d [1 - \mathbb{P}(\tilde{U}_{1s} = 0)]\tilde{\lambda}_s.$$
  \item[b.)] The total claim amount within $i$-th group, $i = 1, 2, ...$ is
$$Y_{i} \stackrel{\rm d}{=} I\{U_{i1} > 0\} \sum_{j = 1}^{U_{i1}} \tilde{Y}_{ij1} + ... + I\{U_{id} > 0\} \sum_{j = 1}^{U_{id}} \tilde{Y}_{ijd},$$
and these r.vs. are strictly positive, i.i.d., having common non-lattice c.d.f., and if $\mathbb{E}\tilde{Y}_{11s} < \infty$ and $\mathbb{E} \tilde{U}_{1s} < \infty$, $d = 1, 2, ..., d$ their mean is
$$\mathbb{E}(Y_1) =  \sum_{s = 1}^d\frac{\tilde{\lambda}_s\mathbb{E}\tilde{U}_{1s}}{(\tilde{\lambda}_1 + ... + \tilde{\lambda}_d)(1 - p_0)} \mathbb{E} \tilde{Y}_{11s} < \infty.$$
In case, when the coordinates of $\overrightarrow{\tilde{Y}}_{111}$ are independent
$$\mathbb{E}(e^{zY_1}) = \sum_{s = 1}^d \frac{\tilde{\lambda}_sg_{\tilde{U}_{1s}}(\mathbb{E}(e^{z\tilde{Y}_{1s}}))-\mathbb{P}(\tilde{U}_{1s} = 0)}{(\tilde{\lambda}_1 + \tilde{\lambda}_2 + ... + \tilde{\lambda}_d)(1 - p_0)}.$$
%$$\mathbb{E}(Y_1^2) = \sum_{s = 1}^d \frac{\tilde{\lambda}_s\left\{\mathbb{D} \tilde{U}_{1s}[\mathbb{E} \tilde{Y}_{1s}]^2 + (\mathbb{E} \tilde{U}_{1s})^2[\mathbb{E} \tilde{Y}_{1s}]^2+ \mathbb{E} \tilde{U}_{1s}\mathbb{D} \tilde{Y}_{1s}\right\}}{(\tilde{\lambda}_1 + \tilde{\lambda}_2 + ... + \tilde{\lambda}_d)(1 - p_0)},$$
If $\mathbb{D}\tilde{Y}_{11s} < \infty$ and $\mathbb{D} \tilde{U}_{1s} < \infty$, $d = 1, 2, ..., d$, then
\begin{eqnarray*}
  \mathbb{D}(Y_1) &=& \sum_{s = 1}^d \frac{\tilde{\lambda}_s\left\{\mathbb{D} \tilde{U}_{1s}[\mathbb{E} \tilde{Y}_{1s}]^2 + (\mathbb{E} \tilde{U}_{1s})^2[\mathbb{E} \tilde{Y}_{1s}]^2+ \mathbb{E} \tilde{U}_{1s}\mathbb{D} \tilde{Y}_{1s}\right\}}{(\tilde{\lambda}_1 + \tilde{\lambda}_2 + ... + \tilde{\lambda}_d)(1 - p_0)} \\
   &-& \left\{\sum_{s = 1}^d\frac{\tilde{\lambda}_s\mathbb{E}\tilde{U}_{1s}}{(\tilde{\lambda}_1 + ... + \tilde{\lambda}_d)(1 - p_0)} \mathbb{E} \tilde{Y}_{11s}\right\}^2.
\end{eqnarray*}
In case when $d = 1$
$$\mathbb{D}(Y_1) =  \frac{\mathbb{D} \tilde{U}_{11}[\mathbb{E} \tilde{Y}_{11}]^2 + \mathbb{E} \tilde{U}_{11}\mathbb{D} \tilde{Y}_{11}}{1 - p_0} -\frac{p_0(\mathbb{E}\tilde{U}_{11})^2(\mathbb{E} \tilde{Y}_{11})^2}{(1 - p_0)^2}.$$
%$$FI\,\,Y_1 = \left[FI\,\,\tilde{U}_{11} - \frac{p_0\mathbb{E}\tilde{U}_{11}}{1 - p_0}\right]\mathbb{E} \tilde{Y}_{11} + FI\,\, \tilde{Y}_{11}.$$

Denote the c.d.f. of $\tilde{Y}_{11s}$ with $F_s$. Then for $x \geq 0$,
$$F_{Y_1}(x) =   \sum_{s = 1}^d \sum_{u = \tilde{L}_s}^{\tilde{R}_s} F_s^{u*}(x) \frac{\tilde{\lambda}_s\mathbb{P}(\tilde{U}_{1s} = u)}{(\tilde{\lambda}_1 + ... + \tilde{\lambda}_d)(1-p_0)},$$
where $0 \leq \tilde{L}_s$ and $\tilde{R}_s = R_s$ are correspondingly the lower and the upper end points of the distribution of
$\tilde{U}_{1s}$, $s = 1, 2, ..., d$.
\begin{equation}\label{cdfY1GroupIndependent}
\bar{F}_{Y_1}(x) =\sum_{s = 1}^d \sum_{u = \tilde{L}_s}^{\tilde{R}_s} [1- F_s^{u*}(x)] \frac{\tilde{\lambda}_s\mathbb{P}(\tilde{U}_{1s} = u)}{(\tilde{\lambda}_1 + ... + \tilde{\lambda}_d)(1-p_0)}.
\end{equation}

The independence between the coordinates is important for the following form of the Laplace - Stieltjes transform
\begin{eqnarray}
  l_{Y_1}(z) &=& g_{U_{11}, ..., U_{1d}}(l_{\tilde{Y}_{111}}(z), l_{\tilde{Y}_{112}}(z), ..., l_{\tilde{Y}_{11d}}(z)) \\\label{LStY1GRoupInd}
\nonumber   &=& \frac{\sum_{s=1}^d \frac{\tilde{\lambda}_s}{\tilde{\lambda}_1 + \tilde{\lambda}_2 + ... + \tilde{\lambda}_d}g_{\tilde{U}_{1s}}(l_{\tilde{Y}_{11s}}(z))- p_0}{1 - p_0}
\end{eqnarray}
and
$$l_{F_{I, Y_1}}(z) = \frac{\sum_{s=1}^d \tilde{\lambda}_s[1-g_{\tilde{U}_{1s}}(l_{\tilde{Y}_{11s}}(z))]}{z\sum_{s = 1}^d\tilde{\lambda}_s\mathbb{E}\tilde{U}_{1s}\mathbb{E} \tilde{Y}_{11s}}.$$
 \end{description}
 \vspace*{12pt}

{\bf Corollary~5.1}. For the risk process defined in (\ref{GeneralModelGroupsIndependent}) if $\mathbb{E}\tilde{Y}_{11s} < \infty$ and $\mathbb{E} U_{1s} < \infty$, $d = 1, 2, ..., d$, then
\begin{description}
 \item[a.)] the safety loading  is  $\rho = \frac{c}{\sum_{s=1}^d\tilde{\lambda}_s\mathbb{E}\tilde{U}_{1s}\mathbb{E}\tilde{Y}_{11s}} - 1$ and net profit condition is $\frac{c}{\sum_{s=1}^d\tilde{\lambda}_s\mathbb{E}\tilde{U}_{1s}\mathbb{E}\tilde{Y}_{11s}} > 1.$
\item[b.)] $\psi(0) = \frac{\sum_{s=1}^d\tilde{\lambda}_s\mathbb{E}\tilde{U}_{1s}\mathbb{E}\tilde{Y}_{11s}}{c}$ and $\delta(0) = 1 - \frac{\sum_{s=1}^d\tilde{\lambda}_s\mathbb{E}\tilde{U}_{1s}\mathbb{E}\tilde{Y}_{11s}}{c}$.
\item[c.)] In case when the coordinates of $\overrightarrow{\tilde{Y}}_{111}$ are independent,   the distribution  of the deficit at the time of ruin with initial capital $0$ is given by

    $\mathbb{P}(-R(\tau(0)+)\leq x| \tau(0) < \infty) = $
$$ = \frac{1}{\sum_{s = 1}^d\tilde{\lambda}_s\mathbb{E}\tilde{U}_{1s}\mathbb{E} \tilde{Y}_{11s}}\sum_{s = 1}^d \sum_{j = \tilde{L}_s}^{\tilde{R}_s}  \tilde{\lambda}_s\mathbb{P}(\tilde{U}_{1s} = j)\int_0^x [1- F_s^{j*}(y)]dy.$$
and if additionally $\mathbb{D}\tilde{Y}_{1s} < \infty$ and $\mathbb{D} U_{1s} < \infty$, $d = 1, 2, ..., d$

$\mathbb{E}(-R(\tau(0)+)| \tau(0) < \infty) =$
$$ = \frac{\sum_{s = 1}^d \tilde{\lambda}_s\left\{\mathbb{D} \tilde{U}_{1s}[\mathbb{E} \tilde{Y}_{1s}]^2 + (\mathbb{E} \tilde{U}_{1s})^2[\mathbb{E} \tilde{Y}_{1s}]^2+ \mathbb{E} \tilde{U}_{1s}\mathbb{D} \tilde{Y}_{1s}\right\}}{2\sum_{s = 1}^d \tilde{\lambda_s}\mathbb{E}\tilde{U}_{1s}\mathbb{E} \tilde{Y}_{1s}}.$$
And for $d = 1$ $$ \mathbb{E}(-R(\tau(0)+)| \tau(0) < \infty) = \frac{1}{2} \left[(FI \tilde{U}_{11} + \mathbb{E}\tilde{U}_{11})\mathbb{E} \tilde{Y}_{1s} + FI \tilde{Y}_{11}\right].$$
\item[d.)] For $\mathbb{D}\tilde{Y}_{1s} < \infty$ and $\mathbb{D} U_{1s} < \infty$, $d = 1, 2, ..., d$ the equation (\ref{tau0}) gives us immediately

$\mathbb{E}(\tau(0)|\tau(0) < \infty) = $
$$ = \frac{\sum_{s = 1}^d \tilde{\lambda}_s\left\{\mathbb{D} \tilde{U}_{1s}[\mathbb{E} \tilde{Y}_{1s}]^2 + (\mathbb{E} \tilde{U}_{1s})^2[\mathbb{E} \tilde{Y}_{1s}]^2+ \mathbb{E} \tilde{U}_{1s}\mathbb{D} \tilde{Y}_{1s}\right\}}{2\left(\sum_{s = 1}^d \tilde{\lambda_s}\mathbb{E}\tilde{U}_{1s}\mathbb{E} \tilde{Y}_{1s}\right)\left(c-\sum_{s = 1}^d \tilde{\lambda_s}\mathbb{E}\tilde{U}_{1s}\mathbb{E} \tilde{Y}_{1s}\right)}.$$
And for $d = 1$
 $$ \mathbb{E}(\tau(0)|\tau(0) < \infty) = \frac{(FI\tilde{U}_{11} + \mathbb{E}\tilde{U}_{11})\mathbb{E} \tilde{Y}_{11} + FI \tilde{Y}_{11}}{2\left(c - \tilde{\lambda}\mathbb{E}\tilde{U}_{11} \mathbb{E} \tilde{Y}_{11}\right)} .$$
 \item[e.)]The joint distribution of the severity of (deficit at) ruin and the risk surplus just before the ruin with initial capital zero is:

 $\mathbb{P}(-R(\tau(0)+) > x, R(\tau(0)-) > y|\tau(0) < \infty) =$
$$ = \frac{\sum_{s = 1}^d \sum_{j = \tilde{L}_s}^{\tilde{R}_s}  \tilde{\lambda}_s\mathbb{P}(\tilde{U}_{1s} = j)\int_{x+y}^\infty [1- F_s^{j*}(z)]dz}{\sum_{s = 1}^d\tilde{\lambda}_s\mathbb{E}\tilde{U}_{1s}\mathbb{E} \tilde{Y}_{11s}}.$$
\item[f.)] The distribution of the claim causing ruin (\ref{theclaimcausingruin}) in this case is

$\mathbb{P}(R(\tau(0)-)-R(\tau(0)+) \leq x| \tau(0) < \infty) =$
$$ = \frac{1}{\sum_{s = 1}^d\tilde{\lambda}_s\mathbb{E}\tilde{U}_{1s}\mathbb{E} \tilde{Y}_{11s}}\sum_{s = 1}^d \sum_{j = \tilde{L}_s}^{\tilde{R}_s}\tilde{\lambda}_s\mathbb{P}(\tilde{U}_{1s} = j)\int_0^x y d F_s^{j*}(y).$$
    \item[g.)] For $u > 0$, $\delta(u)$ satisfy the following equation
 \begin{equation}\label{deltaPrimGroupsIndependent}
 \delta'(u) = \frac{(\tilde{\lambda}_1 + \tilde{\lambda}_2 + ... + \tilde{\lambda}_d)(1 - p_0)}{c}.
 \end{equation}
 \end{description}
 $$.\left[\delta(u) - \sum_{s = 1}^{d}  \frac{\tilde{\lambda}_s}{(\tilde{\lambda}_1 + ... + \tilde{\lambda}_d)(1-p_0)}   \sum_{j = \tilde{L}_s}^{R_s}\mathbb{P}(\tilde{U}_{1s} = j)\int_0^u \delta(u - y) d F_s^{j*}(y)\right].$$
  \begin{description}
     \item[h.)] In case when the coordinates of $\overrightarrow{\tilde{Y}}_{111}$ are independent, the Laplace-Stieltjes transform of $\delta(u)$ is
   $$l_{\delta}(s) = \frac{1 - \frac{\sum_{s=1}^d\tilde{\lambda}_s\mathbb{E}\tilde{U}_{1s}\mathbb{E}\tilde{Y}_{11s}}{c}}{1 - \frac{1}{cs} \sum_{s=1}^d \tilde{\lambda}_s[1-g_{\tilde{U}_{1s}}(l_{\tilde{Y}_{11s}}(z))]}.$$
 \item[i.)] For all $u \geq 0$, the defective renewal equation is
\begin{eqnarray*}
 \delta(u)& =& 1 - \frac{1}{c}\sum_{s=1}^d\tilde{\lambda}_s\mathbb{E}\tilde{U}_{1s}\mathbb{E}\tilde{Y}_{11s} + \\
 \nonumber &+& \frac{1}{c}\sum_{s = 1}^d \sum_{j = \tilde{L}_s}^{\tilde{R}_s} \tilde{\lambda}_s\mathbb{P}(\tilde{U}_{1s} = j)\int_0^u \delta(u-y)[1- F_s^{j*}(y)] dy.
\end{eqnarray*}
and
\begin{eqnarray*}
 \nonumber \psi(u) &=& \frac{1}{c}\sum_{s=1}^d\tilde{\lambda}_s\mathbb{E}\tilde{U}_{1s}\mathbb{E}\tilde{Y}_{11s}-\\
 \nonumber  &-&\frac{1}{c}\sum_{s = 1}^d \sum_{j = \tilde{L}_s}^{\tilde{R}_s}  \tilde{\lambda}_s\mathbb{P}(\tilde{U}_{1s} = j)\int_0^x [1- F_s^{j*}(y)]dy \\
 &+&  \frac{1}{c}\sum_{s = 1}^d \sum_{j = \tilde{L}_s}^{\tilde{R}_s} \tilde{\lambda}_s\mathbb{P}(\tilde{U}_{1s} = j)\int_0^u \psi(u-y)[1- F_s^{j*}(y)] dy.
\end{eqnarray*}
\item[k.)]
$$G(0, y) = c\sum_{s = 1}^d \sum_{j = \tilde{L}_s}^{\tilde{R}_s}  \tilde{\lambda}_s\mathbb{P}(\tilde{U}_{1s} = j)\int_0^x [1- F_s^{j*}(y)]dy.$$
Moreover
$$G(u, y) = \frac{1}{c}\left[\sum_{s = 1}^d \sum_{j = \tilde{L}_s}^{\tilde{R}_s}  \tilde{\lambda}_s\mathbb{P}(\tilde{U}_{1s} = j)\int_0^u G(u-x,y) [1- F_s^{j*}(x)]dx +\right.$$
$$+ \left.\sum_{s = 1}^d \sum_{j = \tilde{L}_s}^{\tilde{R}_s}  \tilde{\lambda}_s\mathbb{P}(\tilde{U}_{1s} = j)\int_u^{u+y}[1- F_s^{j*}(x)]dx \right], \,\, u \geq 0.$$
\item[l.)] The small claim condition  means that these exists the Cramer-Lundberg exponent $\epsilon > 0$:
$$\sum_{s = 1}^d \tilde{\lambda}_s[g_{\tilde{U}_{1s}}(l_{\tilde{Y}_{11s}}(-\epsilon)) - \mathbb{P}(\tilde{U}_{1s} = 0)] = c\epsilon.$$
and for this $\epsilon$ the inequality (\ref{psiLundbergexponent}) is satisfied, i.e. for all $u\geq 0$, $\psi(u) \leq e^{-\epsilon u}$.
\item[m.)]  If the small claim condition is satisfied and
$$ \frac{\sum_{s = 1}^d \sum_{j = \tilde{L}_s}^{\tilde{R}_s}\tilde{\lambda}_s\mathbb{P}(\tilde{U}_{1s} = j)\int_0^\infty x e^{\epsilon x}[1 - F_s^{j*}(x)]dx}{\sum_{s = 1}^d \tilde{\lambda}_s\mathbb{E}\tilde{U}_{1s} \mathbb{E}\tilde{Y}_{11s} }< \infty,$$
then the following Cramer-Lundberg approximation of the probability of ruin holds
$$
\lim_{u \to \infty} e^{\epsilon u}\psi(u) = \frac{c - \sum_{s = 1}^d \tilde{\lambda}_s \mathbb{E}\tilde{U}_{1s} \mathbb{E} \tilde{Y}_{1s}}{\epsilon\sum_{s = 1}^d \sum_{j = \tilde{L}_s}^{\tilde{R}_s}\tilde{\lambda}_s\mathbb{P}(\tilde{U}_{1s} = j)\int_0^\infty x e^{\epsilon x}[1 - F_s^{j*}(x)]dx}.
$$
\item[n.)] If $Y_1$ belongs to the class $D$ of dominatedly varying distributions, given the net profit condition is satisfied,
$$\lim_{u \to \infty} \frac{\psi(u)}{\sum_{s = 1}^d \sum_{j = \tilde{L}_s}^{\tilde{R}_s}  \tilde{\lambda}_s\mathbb{P}(\tilde{U}_{1s} = j)\int_{u}^\infty [1- F_s^{j*}(z)]dz} =$$
$$ = \frac{\sum_{s = 1}^d\tilde{\lambda}_s\mathbb{E}\tilde{U}_{1s}\mathbb{E} \tilde{Y}_{11s}}{c-\sum_{s = 1}^d\tilde{\lambda}_s\mathbb{E}\tilde{U}_{1s}\mathbb{E} \tilde{Y}_{11s}}.$$
\end{description}
  \vspace*{12pt}

\section{Some particular cases}

In this section we describe some examples of recently developed risk models which are particular cases of the models discussed here and therefore particular cases of the Cramer-Lundberg model.

{\bf Example 1.}  The Poisson risk process of order k. It is partially investigated by Kostadinova \cite{kostadinova2013poisson}. In this model $d = 1$ and $U_1$ and $\tilde{U}_1$ are Uniformly distributed over the numbers $1, 2, ..., k$. The last means that we have no possibility for empty group and the number of claims can be any of the numbers $1, 2, ..., k$ with equal probabilities.

{\bf Example 2.}  Polya-Aeppli or order k risk model investigated  in Chukova and Minkova \cite{chukova2015polya}. In their model $\tilde{U}_{11}, \tilde{U}_{21}, ...$ and $U_{11}, U_{21}, ...$ are discrete truncated geometrically distributed over the points $\{1, 2, ..., k\}$ and independent r.vs.
      $$P(U_{11} = i) = \frac{p}{1 - (1-p)^k}(1-p)^{i-1}, \,\,g_{U_1}(s) = \frac{zp(1 - p)}{1-(1 - p)^k}.\frac{1 - (1-p)^kz^k}{1-z(1-p)}, $$
 $$ \mathbb{E}U_{11} = \frac{1 - (1 - p)^k(1+kp)}{p[1-(1 - p)^k]}.$$
      Again we have no possibility for empty group of claims.

{\bf Example 3.}  "Multivariate"  Poisson negative binomial risk process. For $d = 1$ this model is partially investigated in Kostadinova \cite{KrasiDoctoralConference}. In this model the numbers of the claims within the $1, 2, ..., d$-th groups are dependent and  $(\tilde{U}_{11}$,..., $\tilde{U}_{1d})$, $(\tilde{U}_{21},...$, $\tilde{U}_{2d})$, ... are i.i.d. Negative multinomially distributed with parameters$(n; p_1, p_2, ..., p_d)$, $n \in N$, $p_s \in (0, 1)$, $s = 1, 2, ..., d$ and $p_1 + p_2 + ... + p_d < 1$. This model is a particular case of those considered in Section 4.

{\bf Example 4.}   "Multivariate" Poisson negative binomial risk process with independent coordinates of the counting processes. Here $\tilde{U}_{1s}, \tilde{U}_{2s}, ...$, $s = 1, 2, ..., d$ are negative binomially distributed with parameters $(n_s, p_s)$, $n_s \in N$, $p_s \in (0, 1)$, $s = 1, 2, ..., d$. In this case we have possibility for empty group. This model is a particular case of those considered in Section 5. For $d = 1, 2$ it is particularly investigated in Kostadinova and Minkova \cite{kostadinova2014bivariate}.

{\bf Example 5.} Polya-Aeppli risk model. In these models $\tilde{U}_{1s}, \tilde{U}_{2s}, ...$ , $s = 1, 2, ..., d$ are i.i.d. Shifted negative binomially distributed with parameters $(n_s$, $p_s)$ and the  coordinates of the counting processes are independent. For $d = 1$ and $n = 1$ it is introduced and partially investigated by Minkova \cite{minkova1900polya} and considered again in \cite{minkova2009compound}. For $d = 2$ see Kostadinova \cite{kostadinovapolya}. The theory for Polya-Aeppli distributions could be seen e.g. in Johnson, Kotz and Kemp (1992).

{\bf Example 6.} Compound compound Poisson risk model considered in Minkova  \cite{minkova2009compound}. In that model $d = 1$ and $\tilde{U}_{11}$ has arbitrary discrete distribution on the non-negative integers.

{\bf Example 7.} The Poisson model for combining different lines of business, introduced in Wang \cite{Wang1998} is a particular case of the model discussed in Section 4 because he allows no more than one claim to arrive within a group of a fixed type. Another generalization of the model of Wang is that here we allow empty groups.

{\bf Example 8.} The Poisson model with common shock, considered in Cossete, Marceau \cite{A1} is a particular case of the main model considered here with three types of claims, $d = 3$. The claim arrival process is a generalization of Marshall-Olkin process \cite{MO1967, MO1988}. See Figure \ref{fig:Fig4}. The processes
\begin{eqnarray*}
  N_{c,1}  & \sim & HPP(\lambda_{11} + \lambda_{12} + \lambda_{13} + \lambda_{123})\\
  N_{c,2}  & \sim & HPP(\lambda_{21} + \lambda_{22} + \lambda_{23} + \lambda_{123})\\
 N_{c,3}  & \sim & HPP(\lambda_{31} + \lambda_{32} + \lambda_{33} + \lambda_{123})
\end{eqnarray*}
are dependent.
\begin{center}
\includegraphics[width=.9\textwidth]{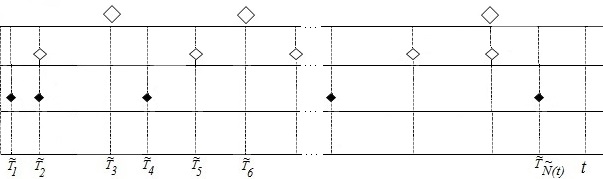}

{\small Fig. 4. Poisson model with common shock. \label{fig:Fig4}}
\end{center}

No more than one claim of a fixed type is possible within a group. There is no possibility the group of all types to be simultaneously empty.  The distribution of the claim numbers is the following
$$P(U_{11} = 1, U_{12} = 0, U_{13} = 0) = \frac{\lambda_{11}}{\lambda_0}, \quad P(U_{11} = 0, U_{12} = 1, U_{13} = 0) = \frac{\lambda_{22}}{\lambda_0},$$
$$P(U_{11} = 0, U_{12} = 0, U_{13} = 1) = \frac{\lambda_{33}}{\lambda_0}, \quad P(U_{11} = 1, U_{12} = 1, U_{13} = 0) = \frac{\lambda_{12}}{\lambda_0},$$
$$P(U_{11} = 1, U_{12} = 0, U_{13} = 1) = \frac{\lambda_{13}}{\lambda_0}, \quad P(U_{11} = 0, U_{12} = 1, U_{13} = 1) = \frac{\lambda_{23}}{\lambda_0},$$
$$P(U_{11} = 1, U_{12} = 1, U_{13} = 1) = \frac{\lambda_{123}}{\lambda_0},$$
where $\lambda_0 = \lambda_{11} + \lambda_{22} + \lambda_{33} + \lambda_{13} + \lambda_{12} + \lambda_{23} + \lambda_{123}$ and $N(t) \sim HPP (\lambda_0)$.

\section*{Acknowledgements}
The work was  supported by project Fondecyt Proyecto Regular No. 1151441, Project LIT-2016-1-SEE-023,
and by the bilateral projects Bulgaria - Austria, 2016-2019, Contract number 01/8, 23/08/2017.

%\newpage

\end{document}